\documentclass[brochure,12pt]{bourbaki}
\usepackage[matrix,arrow]{xy}
\usepackage{amssymb,amsfonts,amsmath,footnote}
\usepackage[francais]{babel}
\usepackage{graphics}
\usepackage{graphicx}
\usepackage{psfrag}
\usepackage{color}
\newcommand{\be}{\begin{equation}}
\newcommand{\ee}{\end{equation}}
\newcommand{\R}{\mathbb{R}}
\newcommand{\N}{\mathbb{N}}
\newcommand{\C}{\mathbb{C}}
\newcommand{\Z}{\mathbb{Z}}

\newcommand\res{\mathop{\hbox{\vrule height 7pt width .5pt depth 0pt
\vrule height .5pt width 6pt depth 0pt}}\nolimits}
\def\thesection{\Roman{section}}
\def\theequation{\thesection.\arabic{equation}}

\def\ti{\tilde}
\def\lf{\left}
\def\rg{\right}

\def\al{\alpha}
\def\la{\lambda}

\def\ep{\varepsilon}
\def\ds{\displaystyle}
\def\ov{\overline}

\def\om{\omega}
\def\p{\partial}

\def\res{\mathop{\hbox{\vrule height 7pt width .5pt 
depth 0pt\vrule height .5pt width 6pt depth 0pt}}\nolimits}

\date{Janvier 2014}
\bbkannee{66\`eme ann\'ee, 2013-2014}
\bbknumero{1081}
\title{M\'ethodes de Min-Max et la Conjecture de Willmore}
\subtitle{d'apr\`es F.C.Marques et A.A.Neves }
\author{Tristan RIVI\`ERE}
\address{ETH-Zentrum\\
Forschungsinstitut f\"ur Mathematik\\
CH--8092 Z\"urich, Switzerland}
\email{tristan.riviere@math.ethz.ch}

\begin{document}
\maketitle

\noindent{\bf INTRODUCTION}

\bigskip
L'\'etude des variations de Lagrangiens, c'est \`a dire la recherche de leurs points critiques, de leur indice, de la topologie de leurs ensembles de niveau...etc sont des probl\'ematiques anciennes des math\'ematiques qui vont bien au-del\`a du stricte champ du calcul des variations et de l'analyse en g\'en\'eral. Elles ont eu des impacts importants dans bien d'autres domaines  comme la topologie diff\'erentielle, la g\'eom\'etrie diff\'erentielle et riemannienne ou encore la g\'eom\'etrie alg\'ebrique complexe...etc (th\'eorie de Morse, th\'eorie de Jauges, espaces de modules..etc).

Une probl\'ematique d\'elicate du calcul des variations consiste \`a identifier des points critiques autrement que
par simple proc\'edure de minimisation, en d'autre termes, la recherche de points selles. 
Au d\'ebut du XXeme sicle G.D. Birkhoff, dans son \'etude de l'existence de g\'eod\'esiques ferm\'ees dans l'espace
des courbes homotopiques \`a un point sur une vari\'et\'e compacte sans bord, a d\'evelopp\'e avec succ\`es une technique dite de ''min-max''. Cette approche telle qu'il l'a conue s'est av\'er\'ee trs performante aussi longtemps que
l'on travaille avec des objets 1-dimensionnels.

La th\'eorie de la mesure g\'eom\'etrique, telle qu'elle a \'emerg\'e au d\'ebut des annŽes 50, fusionnant la th\'eorie des courants de de Rham avec la notion d'ensemble rectifiable \'etudi\'ee en particulier par A.Besicovitsch les d\'ec\'enies prŽcŽdentes, \'etait motiv\'ee
par des probl\`emes g\'en\'eraux de minimisation de volume pour des sous vari\'et\'es (ou plus exactement de leur version affaiblie : des courants entiers rectifiables) sous diff\'erentes contraintes de bord ou plus g\'en\'eralement d'homologie. Ce genre de questions avaient elle m\^eme leur origine dans le fameux probl\`eme de Plateau 
consistant \`a chercher un disque ''optimal'' bordant une courbe de Jordan dans un espace euclidien donn\'e.
La th\'eorie de la mesure g\'eom\'etrique, au del\`a des questions de minimisation, s'est  progressivement  empar\'ee de la question de trouver des points selles et c'est autour des annŽes 70 que F.Almgren et J.Pitts donn\`erent le cadre d'une m\'ethode de min-max, inspir\'ee de Birkhoff, pour des objets tr\`es g\'en\'eraux comme les courants rectifiables.

Il y a bient\^ot 2 ans F.C.Marques et A. Neves ont mis en oeuvre cette m\'ethode de min-max dans le cadre des courants rectifiables ferm\'es de dimension
2 dans la sph\`ere 3-dimensionnelle pour des d\'eformations correspondants, entre autres, \`a l'action du groupe des transformations conformes.
Ils sont ainsi parvenus \`a mettre en \'evidence une surface minimale, de courbure moyenne nulle, qui minimise l'aire parmi toutes
les autres surfaces minimales ferm\'ees de genre non nul. Le calcul pr\'ecis de l'indice de cette surface leur a permis d'identifier cette surface comme \'etant le fameux tore de Clifford. Une des cons\'equences spectaculaire de ce r\'esultat est la d\'emonstration de la conjecture dite ''de Willmore'' qui pr\'edisait la forme optimale des surfaces de genre non nul minimisant le lagrangien du m\^eme nom, introduit dans le cadre de la g\'eom\'etrie conforme par W.Blaschke au d\'ebut du XXeme si\`ecle dans son effort de fusionner l'invariance conforme et la th\'eorie des surfaces minimales.

\section{L'origine des m\'ethodes de Min-Max et la recherche de g\'eod\'esiques f\'erm\'ees.}

Trouver une courbe dans une vari\'et\'e riemannienne $(M,g)$ connexe ferm\'ee minimisant la longueur sous une contrainte 
homologique ou homotopique simple comme par exemple celle de prescrire ses extr\'emit\'es lorsqu'elle n'est pas ferm\'ee
ou comme celle d'appartenir \`a  une classe d'homotopie donn\'ee du $\pi_1(M)$ lorsque la courbe est un cercle est peut \^etre l'exercice le plus 
\'el\'ementaire du calcul des variations. Cet exercice consiste essentiellement \`a s'assurer que l'espace des courbes dans lequel
on pose le probl\`eme est bien ferm\'e pour cette op\'eration de minimisation et \`a faire un petit travail de preuve de r\'egularit\'e
sur le minimum obtenu . 

\medskip

La question devient plus d'\'elicate lorsque l'on se pose la question de trouver une g\'eod\'esique ferm\'ee 
sur un espace simplement connexe comme la sph\`ere $S^2$ \'equip\'ee d'une m\'etrique quelconque. Cette question fut probablement 
consid\'er\'ee pour la premi\`ere fois dans les travaux d'Hadamard \cite{Ha} et de Poincar\'e \cite{Po} vers le d\'ebut du XXeme si\`ecle.

\medskip

G.D.Birkhoff a d\'evelopp\'e une m\'ethode dite de {\it ''min-max''}  afin de d\'emontrer le r\'esultat suivant.

\begin{theo} [\cite{Bi}]
\label{birkhoff}
Toute vari\'et\'e riemannienne hom\'eomorphe \`a la sph\`ere poss\`ede une g\'eod\'esique ferm\'ee non triviale.
\end{theo}

L'id\'ee principale de la preuve est la suivante. On consid\`ere l'espace ${\mathcal C}$ des applications continues de $S^1$ dans $S^2$ r\'eguli\`eres par morceaux que l'on 
quotiente par l'action des hom\'eomorphismes positifs et  r\'eguliers par morceaux de $S^1$, espace que l'on ''compactifie'' en rajoutant les points de $S^2$ vus comme des courbes d\'eg\'en\'er\'ees. Dans
cet espace on consid\`ere le 1-cycle non trivial de $H_1({\mathcal C},S^2,{\Z})$ donn\'e par 
\[
t\in[0,1]\longrightarrow \gamma(t):=S^2\cap\{z=1-2\,t\}
\]
Tout chemin homotope \`a $\gamma(t)$, pour des d\'eformations maintenants les extr\'emit\'es dans les ''courbes d\'eg\'en\'er\'ees'' donn\'ees par les points de $S^2$, constitue ce que l'on appelle un ''{\it balayage}'' de $S^2$ dans le sens qu'il 
r\'ealise un g\'en\'erateur de $H_2(S^2)$. La g\'eod\'esique recherch\'ee sera alors un point critique du lagrangien de longueur $L$ donn\'e par la
proc\'edure de {\it min-max} suivante
\[
c=\inf_{\sigma\simeq\gamma}\ \max_{t\in [0,1]} L(\sigma(t))
\]
o\`u $\simeq$ est l'\'equivalence d'homotopie mentionn\'ee plus haut. La valeur $c$ est appel\'ee la largeur du {\it min max}. On d\'emontre que l'infimum $c$ est atteind grace \`a un argument de comparaison d\'evelopp\'e par Birkhoff et connu sous le nom de {\it ''d\'eformation polynomiale''}
ou {\it ''racourcicement de longueur''} de Birkhoff qui d\'eforme tout \'el\'ement $\sigma$ en un \'el\'ement homotope et r\'ealisant pour tout $t$ des g\'eod\'esiques
par morceau tout en d\'ecroissant le maximum des longeurs (voir \cite{Cr} et \cite{CM} par exemple).

\medskip

En 1929 L.Lyusternik et L.Schnirel'man annonc\`erent dans une note publi\'ee au {\it Comptes Rendus de l'Acad\'emie des Sciences de Paris}, voir \cite{LS}, un r\'esultat
qui d\'emontrait un probl\`eme laiss\'e ouvert par H.Poincar\'e.

\medskip

\begin{theo} 
\label{Lyusternik-Shnirel'man}
Toute vari\'et\'e riemannienne hom\'eomorphe \`a la sph\`ere poss\`ede au moins trois g\'eod\'esiques ferm\'ees non triviales, distinctes les unes des autres,
et ne poss\'edant aucun point d'auto-intersection.
\end{theo}

\medskip

Le r\'esultat est optimal au sens o\`u sur certains ellipsoides il existe exactement trois g\'eod\'esiques ferm\'es distinctes sans point d'auto-intersection (voir \cite{Kl}).
L'article contenant la preuve de ce r\'esultat et publi\'e un an plus tard que la note a ouvert la voie \`a un domaine entier de l'analyse portant sur
les m\'ethodes topologiques pour le calcul des variations. La version anglaise de cet article est paru dans \cite{LS1}. N\'eanmoins la preuve du th\'eor\`eme pr\'ec\'edent n'y \'etait pas compl\`ete et les contributions pour
\'etablir rigoureusement ce r\'esultat ont \'et\'e nombreuses et se sont \'etall\'ees sur plus d'un demi-si\`ecle.
Nous renvoyons le lecteur qui d\'esire en savoir plus sur l'histoire riche en rebondissements de la preuve de ce th\'eor\`eme
\`a l'excellent expos\'e de I.Taimanov (\cite{Ta}).

L'apport principal de l'article de  Lyusternik et Schnirel'man a \'et\'e d'\'etendre la strat\'egie de {\it min-max} de Birkhoff \`a des cycles de dimensions sup\'erieures \`a $H_1({\mathcal C}, S^2)$
et d'introduire des classes d'homotopie non triviales - qu'ils on nomm\'e cat\'egories - allant de $[0,1]^n$ pour $n=1,2,3$ dans l'espace ${\frak C}$ des courbes non-orient\'ees de jordan  rectifiables (c'est \`a dire des images continues et injectives de $S^1$ de
mesure 1-dimensionnelle finie dans $S^2$)  auxquelles on adjoint les courbes d\'eg\'en\'er\'es donn\'ees par les points de $S^2$. Introduire des cycles de dimension sup\'erieure correspond plus ou moins \`a chercher des g\'eod\'esiques d'indice de plus en plus \'elev\'e. 
\[
c_n=\inf_{\sigma\simeq\gamma_n}\ \max_{x\in [0,1]^n} L(\sigma(x))
\]
le chemin $\gamma_1(x_1)$ est le chemin de Birkhoff d\'efini plus haut tandis que $\gamma_2(x_1,x_2)$ et $\gamma_3(x_1,x_2,x_3)$ sont des chemins construits \`a partir de $\gamma_1$ en rajoutant comme
variable suppl\'ementaire l'action des rotations d'axe $x$ et d'angle allant de 0 \`a $\pi$ pour $\gamma_2$ 
\[
\gamma_2(x_1,x_2)=
\left(
\begin{array}{ccc}
1& 0 &0\\[2mm]
0&\cos(\pi x_2) &\sin(\pi x_2)\\[2mm]
0&-\sin(\pi x_2) &\cos(\pi x_2)
\end{array}
\right)
\gamma_1(x_1)
\]
et l'action des rotations respectivement d'axe $x$ et d'axe $y$ toutes deux d'angles allant de 0
\`a $\pi$ pour $\gamma_3$
\[
\gamma_2(x_1,x_2,x_3)=
\left(
\begin{array}{ccc}
\cos(\pi x_3)& 0 &\sin(\pi x_3)\\[2mm]
0&1 &0\\[2mm]
-\sin(\pi x_3)& 0 &\cos(\pi x_3)
\end{array}
\right)
\gamma_2(x_1,x_2)
\]

Par ailleurs la classe des applications $\sigma$ consid\'er\'ees, ainsi que les 3 classes de d\'eformations sur lesquelles les \'equivalence d'homotopies $\simeq$ sont
d\'efinies, sont les 3 espaces d'applications continues de $[0,1]^n$ dans ${\frak C}$ ayant les m\^emes sym\'etries
que $\gamma_n$ respectivement pour $n=1,2,3$ au bord : pour $x_1,x_2\in [0,1/2]$ et $x_3\in [0,1]$
\be
\label{1}
\left\{
\begin{array}{l}
\sigma(x_1,0)=\sigma(1-x_1,1)\quad\mbox{ pour }n=2\\[5mm]
\sigma(x_1,0,x_3)=\sigma(1-x_1, 1,x_3)\\[5mm]
\sigma(x_1,x_2,1)=\sigma(1-x_1,1-x_2,0)
\end{array}
\rg.
\ee
et par ailleurs on impose aussi $\sigma(0)=\gamma_1(0)$ et $\sigma(1)=\gamma_1(1)$ pour $n=1$, $\sigma(0,x_2)=\gamma_2(0,x_2)$  et $\sigma(1,x_2)=\gamma_2(1,x_2)$ pour $n=2$ ainsi que $\sigma(0,x_2,x_3)=\gamma_3(0,x_2,x_3)$  et $\sigma(1,x_2,x_3)=\gamma_3(1,x_2,x_3)$ pour $n=3$. Il convient d'observer que l'espace ${\frak C}$ choisi ne prend
pas en compte l'orientation des courbes ce qui justifie les identit\'es (\ref{1}) en particulier pour $\sigma=\gamma_n$.
Ces trois probl\`emes de {\it min-max} \'etant pos\'es, une grande partie de la difficult\'e de la preuve demeure en cela qu'il s'agit de d\'emontrer que chacun des trois est bien atteint par des g\'eod\'esiques sans auto-intersections. La  technique de {\it ''racourcicement de longueur''} de Birkhoff malheureusement est inutilisable telle qu'elle, celle-ci ne respectant pas la propri\'et\'e
de non auto-intersection. De nombreux travaux, incluant celui de Lyusternik et Schnirel'man, se sont heurt\'es \`a la difficult\'e de trouver un substitut ad-hoc \`a la technique de {\it ''racourcicement de longueur''} de Birkhoff sans parvenir \`a la r\'esoudre compl\`etement dans le cadre g\'en\'eral. L'article \cite{Ta1} a mis fin \`a cette \'errance proposant une preuve compl\`ete de l'atteignabilit\'e de chaqu'un des trois {\it min-max} par des g\'eod\'esiques plong\'ees. Finalement, il n'est pas difficile de d\'emontrer l'in\'egalit\'e
\[
c_1\le c_2\le c_3
\]
Si jamais deux niveaux $c_n$ et $c_{n+1}$ sont \'egaux on d\'emontre alors qu'il existe en fait une infinit\'e de solutions.

\medskip

On vient de voir que pour passer du premier th\'eor\`eme de Birkhoff \`a celui, plus fort sur l'existence de 3 g\'eod\'esiques plong\'ees, afin d'appliquer de proc\'edures de 
{\it min-max} plus \'elabor\'ees, il a fallu affaiblir la notion de courbes orient\'ees r\'eguli\`eres par morceaux \`a celle beaucoup plus g\'en\'erale de cercles rectifiables non orient\'es.

\medskip

L'affaiblissement la notion de courbe pour les probl\`emes de {\it min-max} peut  s'av\'erer\^etre tr\`es efficace comme le montre la d\'emonstration du th\'eor\`eme suivant du \`a E.Calabi.

\begin{theo} 
\label{Calabi}
Sur une vari\'et\'e hom\'eomorphe \`a la sph\`ere toute g\'eod\'esique ferm\'ee de longueur minimale a au plus un point d'auto-intersection.
\end{theo}

La d\'emonstration de ce th\'eor\`eme, telle qu'elle a \'et\'e originellement \'ecrite dans \cite{CC}, fait usage d'un espace encore plus grand que l'espace ${\mathcal C}$ des immersions de $S^1$
r\'eguli\`eres par morceau  ou encore que l'espace ${\mathfrak C}$ des courbes de jordans rectifiables consi\'er\'es plus haut, il s'agit de l'espace ${\mathcal Z}$ des {\it 1-cycles rectifiables}. Cette notion est un des objets fondamentaux constituants la charpente de la th\'eorie de mesure g\'eom\'etrique. Cette th\'eorie s'est av\'er\'ee essentielle pour donner un cadre variationnel
bien pos\'e aux probl\`emes de {\it min-max} d'objets g\'eom\'etriques de dimension plus grande que 1. Il est donc n\'ec\'essaire que nous nous y arr\'etions un moment afin de pouvoir
rendre compte du travail de F.Marques et A.Neves sur la conjecture de Willmore.

\section{Les m\'ethodes de Min-Max en th\'eorie de mesure g\'eom\'etrique d'apr\`es F.Almgren et J.Pitts.}

Le prob\`eme dit {\it de Plateau} tel qu'il a \'et\'e pos\'e vers la fin du XIXeme si\`ecle  par le physicien belge qui lui a donn\'e son nom consiste \`a trouver un disque immerg\'e d'aire minimale bordant
une courbe de Jordan donn\'ee dans un espace euclidien quelconque. Il fut r\'esolu au d\'ebut des ann\'ees 30 pour la premi\`ere fois ind\'ependamment pas J.Douglas et T.Rad\`o, le premier ayant 
re\c cu la m\'edaille Fields en 1936 - la premi\`ere attribu\'ee conjointement avec le math\'ematicien L.Ahlfors. L'aproche de l'\'epoque \'etait alors  une approche dite {\it param\'etrique} qui consiste \`a minimiser
l'aire sur des applications du disque dans l'espace euclidien envoyant continuement et homeomorphiquement le bord du disque sur la courbe rectifiable. Tr\`es vite cette approche \`a montr\'e
ses limites : existerait-il des surfaces de topologie plus compliqu\'ee et d'aire plus basse encore bordant cette m\^eme courbe ? comment \'etendre le probl\`eme en dimension plus grande et trouver
une sous-vari\'et\'e de volume minimal bordant une sous vari\'et'\'e orient\'e ferm\'ee donn\'ee dans un espace euclidien ?...le cadre param\'etrique \'etait trop restreint pour apporter des r\'eponses satisfaisantes.

Cet ensemble de questions m\'elangeant la g\'eom\'etrie et le calcul des variations a donn\'e naissance \`a la {\it th\'eorie de la mesure g\'eom\'etrique} au d\'ebut des ann\'ees 50 sous l'impulsion
de diff\'erents math\'ematiciens parmi lesquels  H.Federer, W.Fleming, E.De Giorgi, E.R.Reifenberg. L'objet de d\'epart de cette th\'eorie sont les {\it courants de de Rham} de dimension $k$ d'une vari\'et\'e riemannienne $(M^m,g)$ donn\'ee qui sont les distributions
vectorielles agissant sur les $k-$ formes $C^\infty$ \`a support compacte. Cet espace est not\'e ${\mathcal D}_k(M^m)$. Le bord $\p C$ d'un tel courant $C$ est d\'efini par
\[
\forall \al\in C^\infty_0(\wedge^{k-1}M^m)\quad \quad\left<\p C,\al\right>=\left<C,d\al\right>
\]
Les courants sont \'equip\'e d'une topologie {\it faible s\'equentielle} qui dit qu'une suite $C_n$ converge faiblement vers $C_\infty$ si
\be
\label{2}
\forall\om\in C^\infty_0(\wedge^{k}M^m)\quad\quad\left<C_n,\om\right>\longrightarrow \left<C_\infty,\om\right>\quad.
\ee
Parmi ces courants, les {\it courants rectifiables entiers} de dimension $k\in {\N}$, dont l'ensemble est not\'e habituellement ${\mathcal R}_k(M^m)$, sont les courants d'int\'egration sur des ensembles rectifiables orient\'es de dimension $k$ - ensembles de mesure ${\mathcal H}^k$ finie et poss\'edant ${\mathcal H}^k$ presque partout\footnote{${\mathcal H}^k$ est la mesure de Hausdorff $k-$dimensionnelle.} un {\it plan tangent approxim\'e} et le choix mesurable d'une orientation de ce plan\footnote{Les ensembles rectifiables peuvent \^etre ''vus'' comme des versions ''th\'eorie de la mesure'' des sous vari\'et\'e de $M^m$.}-  \'equip\'e d'une fonction multiplicit\'e ${\mathcal H}^k$ int\'egrable sur le {\it porteur rectifiable} et \`a valeur enti\`ere. Un tel courant pourra \^etre not\'e $C=(S,\tau,\theta)$ o\`u $S$ est l'ensemble $k-$dimensionnel rectifiable
dit le {\it porteur} du courant, $\tau$ est l'{\it orientation} du courant, c'est une application ${\mathcal H}^k$ mesurable sur $S$ \`a valeur dans le fibr\'e des $k-$multi-vecteurs simples et unit\'e sur $M^m$ et dont la direction $|\tau|$ coincide avec le plan approxim\'e  ${\mathcal H}^k$ presque partout sur $S$ et $\theta$ est la {\it multiplicit\'e}, c'est \`a dire une fonction ${\mathcal H}^k$ mesurable sur $S$ et \`a valeur enti\`ere. L'action d'un tel courant sur une $k-$forme r\'eguli\`ere \`a support compacte $\om$ est donn\'ee par\footnote{L'op\'eration de restriction des courants et mesures se note $\res$ en th\'eorie de mesure g\'eom\'etrique. On \'etant cette notation \`a l'op\'eration de multiplication par une fonction mesurable en g\'en\'eral, la restriction \`a un sous-ensemble mesurable n'\'etant finallement que la mutiplication par la fonction caract\'eristique de cet ensemble.}
\[
<C,\om>=\int_{M^m}\ <\om,\tau>\ \theta\ d{\mathcal H}^k\res S
\]  Le courant d'int\'egration sur une sous-vari\'et\'e orient\'ee de dimension $k$ de $M^m$ est
un exemple d'un tel {\it courant rectifiable entier}. 
Les {\it courant rectifiables entiers} de dimension $k$ dont les bords sont aussi des {\it courant rectifiables entiers} sont appel\'es {\it courants int\'egraux} de dimension $k$ et cet ensemble
sera not\'e ${\mathcal I}_k(M^m)$. Un sous-espace remarquable de ${\mathcal I}_k(M^m)$ est l'espace des {\it cycles rectifiables entiers} que nous avons mentionn\'e pr\'ec\'edemment
et qui coincide avec les \'el\'ements de ${\mathcal I}_k(M^m)$ de bord nul. Cet espace est not\'e ${\mathcal Z}_k(M^m)$.

La comasse d'une $k-$forme $\om$ de $M^m$ est d\'efinie par
\[
\|\om\|_\ast:=\sup_{x\in M^m}\ \sup_{e_1\cdots e_k\in S_xM^m}|<\om,e_1\wedge\cdots\wedge e_k>|\quad,
\]
o\`u $S_xM^m$ d\'esigne la sph\`ere unit\'e de l'espace tangent $(T_xM^m,g)$. La masse d'un courant $k-$dimensionnel $C$ est donn\'ee par
\[
M(C):=\sup\left\{ \left<C,\om\right>\ ;\ \|\om\|_\ast\le 1\right\}\in{\R}^+\cup\{+\infty\}\quad.
\]
La masse du courant d'int\'egration le long d'une sous-vari\'et\'e de $M^m$ coincide avec le volume de cette sous-vari\'et\'e et peut donc \^etre interpr\'et\'ee comme \'etant
la g\'en\'eralisation naturelle \`a $I_k(M^m)$ de cette notion de volume. Plus g\'en\'eralement, pour un {\it courant rectifiable entier} $C=<S,\tau,\theta>$ donn\'e on a
\[
M(<S,\tau,\theta>)=\int_{M^m}|\theta|\ d{\mathcal H}^k\res S\quad.
\]
Le th\'eor\`eme suivant du \`a H.Federer et W.Fleming est reconnu comme \'etant le th\'eor\`eme fondateur de la th\'eorie de la mesure g\'eom\'etrique.

\begin{theo} [\cite{FF}]
\label{FedererFleming}
Soit $C_n$ une suite de courants int\'egraux de $M^m$ satisfaisant
\[
\limsup_{n\rightarrow +\infty} M(C_n)+M(\p C_n)<+\infty
\]
alors il existe une sous-suite $C_{n'}$ qui converge faiblement vers un courant \underbar{int\'egral} $C_\infty$ et par ailleurs
\[
M(C_\infty)\le\liminf M(C_{n'})\quad.
\]
\end{theo}
Grace \`a ce r\'esultat on obtient l'existence d'un solution au {\it probl\`eme de Plateau} pour tout bord int\'egral de dimension $k-1$ fix\'e dans une vari\'et\'e quelconque ou plus g\'en\'eralement
l'existence d'un courant int\'egral de masse minimale dans toute classe d'homologie donn\'ee. La difficult\'e cependant de d\'eterminer la r\'egularit\'e - au-del\`a de la simple rectifiabilit\'e - 
pour chacun de ces minima est la contre-partie de s'\^etre pac\'e dans un espace si grand et a stimul\'e de nombreux travaux tr\`es ardus dont le point culminant est peut-\^etre le r\'esultat
d'Almgren. Ce r\'esultat  affirme que l'ensemble singulier d'un courant rectifiable $k-$dimensionnel d'aire minimale \`a lint\'erieur - dans des ouverts n'intersectant pas son bord -  est de dimension de Hausdorf au plus $k-2$. Une preuve de ce r\'esultat est l'objet
de la publication monumentale, malheureusement non r\'ef\'er\'ee, et posthume de F.Almgren \cite{Al}. Ce r\'esultat a \'et\'e tout nouvellement d\'emontr\'e, en grande partie ind\'ependemment du travail d'Almgren , dans une s\'erie
d'articles remarquables  de C.De Lellis et E.Spadaro \cite{DS1}, \cite{DS2} et \cite{DS3}.

\medskip

En ce qui concerne la r\'egularit\'e des minima, le cas de la codimension 1, $k=m-1$, est particulier et \underbar{relativement} moins difficile que le cas g\'en\'eral de codimension quelconque.
Cela est du au fait que la diversit\'e des {\it cone tangents singuliers}  - apr\`es dilatation infinie en un point - est beaucoup plus restreinte pour les hypersurfaces d'aire minimale. Un r\'esultat
c\'el\`ebres de J.Simons \cite{Si} affirme que tout {\it courant int\'egral} de codimension 1 et d'aire minimale dans une vari\'et\'e de dimension plus petite que 8 est r\'egulier \`a l'int\'erieur. Ce r\'esultat n'est plus vrai en dimension 8 : E.Bombieri, E. De Giorgi et E. Giusti, \cite{BDG}, produisirent un cone singulier \`a l'origine dans ${\R}^8$ et d'aire minimale.

\medskip

 Vers le milieu des ann\'ees 60, l'ambition de produire des surfaces non absolument minimisantes mais points critiques de la masse
germa naturellement dans le cadre de la th\'eorie de la mesure g\'eom\'etrique. Autant le cadre des courants int\'egraux avait \'et\'e tr\`es efficace dans la recherche de minima absolu
autant il semblait inadapt\'e pour d\'evelopper des m\'ethodes de {\it min-max}. En effet, la masse n'est pas continue, mais seulement semi-continue inf\'erieurement, dans l'espace des courants int\'egraux munis de la topologie faible s\'equentielle. En 1965, F.Almgren, produisit un article tr\`es innovant\footnote{Cet article fondamental dans la th\'eorie, pour des raisons un peu surprenantes  qu'explique B.White dans son article en hommage \`a son ancien directeur de th\`ese disparu \cite{Wh}, n'a jamais \'et\'e publi\'e.}, introduisant une notion nouvelle de ''surfaces faibles'' : les varifolds\footnote{L'auteur de cet expos\'e s'excuse pour ce ''franglais'' mais n'a trouv\'e aucune traduction satisfaisante au concept introduit par F.Almgren. L'id\'ee d'Almgren \'etait de produire une notion ad-hoc de vari\'et\'e - {\it manifold} en anglais - compatible avec le
calcul des variations, d'o\`u la fusion {\it var-ifold}. En fran\c cais une telle fusion donnerait {\it var-i\'et\'e}...ce qui tombe vraiment mal il faut le reconna\^ \i tre.}. Ces objets math\'ematiques avaient en fait d\'ej\`a \'et\'e consid\'er\'es par L.C Young une quinzaine d'ann\'ees auparavant sous le nom de ''surfaces g\'en\'eralis\'ees'' et est tr\`es proche de la notion connue sous le nom de {\it mesure de Young}. F.Almgren cependant a pouss\'e bien plus loin cette notion pour les besoins du calcul des variations.

\medskip

Nous notons $G_k(M^m)$ le fibr\'e en grassmanienne des plans $k-$dimensionnels non orient\'es de $TM$. Un {\it varifold} $k-$dimensionnel dans la vari\'et\'e riemannienne $(M^m,g)$
est une mesure de Radon sur $G_k(M^m)$. L'espace des mesures de Radon sur $G_k(M^m)$ \'equip\'e de la topologie faible $\ast$ en dualit\'e avec les fonctions continues de $G_k(M^m)$
et not\'e $V_k(M^m)$. Par exemple si $S$ est un ensemble rectifiable $k-$dimensionnel de $M^m$ on note $|S|$ l'\'el\'ement de $V_k(M^m)$ donn\'e par
\[
\forall A\ \mbox{ensemble mesurable de }G_k(M^m)\quad |S|(A)={\mathcal H}^k\left(S\cap\left\{ x\ ;\ Tan^k({\mathcal H}^k\res S,x)\in A\right\}\right)
\]
o\`u $Tan^k({\mathcal H}^k\res S,x)$ est la mesure tangente obtenue en dilatant la mesure de Hausdorff $k-$dimensionnelle restreinte \`a $S$ au point $x$ et qui coincide
avec la mesure de Hausdorff restreinte au plan tangent approxim\'e \`a $S$ en $x$. L'hypoth\`ese de rectifiabilit\'e consiste exactement \`a dire qu'une telle limite $Tan^k({\mathcal H}^k\res S,x)$ existe ${\mathcal H}^k$-presque partout sur $S$.

\'Etant donn\'e un {\it varifold} quelconque $V$ de $V_k(M^m)$ on lui associe une mesure de Radon sur $M^m$ not\'ee $\|V\|$ et d\'efinie par
\[
\forall B\ \mbox{ensemble mesurable de }M^m\quad \|V\|(B)=|V|(\pi^{-1}(B))\quad,
\]
o\`u $\pi$ est la projection du fibr\'e en grassmanienne $G_k(M^m)$ sur la base $M^m$. En particulier si $V=|S|$ est le {\it varifold} associ\'e \`a l'ensemble rectifiable $S$
tel que nous venons de le d\'efinir plus haut, on a $\||S|\|={\mathcal H}^k\res S$.

On d\'enote $IV_k(M^m)$ l'espace des {\it varifold entiers rectifiables} de dimension $k$ aussi appel\'e {\it varifolds int\'egraux} de dimension $k$ : c'est \`a dire le sous espace de $V_k(M^m)$ des mesures de Radons de $G_k(M^m)$ obtenues comme une somme convergente de mesures de la forme
$|S|\res\theta$ o\`u $S$ est un sous ensemble $k-$dimensionnel rectifiable de $M^m$ quelconque et $\theta$ est une fonction ${\mathcal H}^k$ mesurable sur $S$ quelconque et \`a valeur dans $\N$. Si $C=<S,\tau,\theta>$ est un {\it courant rectifiable entier} on d\'efinit son {\it varifold associ\'e} par
\[
|C|:=|S|\res |\theta|
\]
La fermeture des {\it varifold entiers rectifiables} de dimension $k$ pour la topologie faible $\ast$ avec les fonctions continues sur $G_k(M)$ est not\'e ${\mathcal V}_k(M)$.

On v\'erifie sans difficult\'es que si $V$ est un {\it varifold entier rectifiable} alors $\|V\|$ poss\`ede $\|V\|$ presque partout une densit\'e\footnote{Grace \`a un th\'eor\`eme difficile de D.Preiss, \cite{Pr} nous savons qu'une mesure Bor\'elienne non n\'egative $\mu$ poss\`ede une densit\'e $k-$dimensionnelle $\mu$ presque partout si est seulement si elle est port\'ee par un ensemble $k-$dimensionnel rectifiable : il existe un ensemble $k-$dimensionnel rectifiable $S$ est une fonction positive et ${\mathcal H}^k$ mesurable $f$ telle que $\mu=f\, {\mathcal H}^k\res S$. Ce r\'esultat r\'epondait
par la positive \`a une conjecture de  A.Besicovitch dont les travaux, principalement entre les deux guerres mondiales, sont \`a l'origine du d\'eveloppement de la th\'eorie de la mesure g\'eom\'etrique. Ce r\'esultat de Preiss a eu des impacts important dans l'\'etude des ensembles singuliers de certaines EDP comme Yang-Mills ou les applications harmoniques.}
\[
\Theta^k(\|V\|,x):=\lim_{r \rightarrow 0}r^{-k}\|V\|(B_r(x))\quad\mbox{ existe pour }\|V\|\mbox{ presque tout }x\in M^m
\]
o\`u $B_r(x)$ est la boule g\'eod\'esique de centre $x$ et de rayon $r$  et $\|V\|$ poss\'ede une mesure tangente \'egale \`a $\theta_0\ {\mathcal H}^k\res P$ \`ou $P$ est un plan tangent $k-$dimensionnel au point consi\'er\'e et $\theta_0\in {\N}^\ast$.
 
L'avantage principal des {\it varifolds} sur les courants est la continuit\'e de la masse pour la topologie faible s\'equentielle lorsque $M^m$ est une vari\'et\'e compacte. On d\'eduit en effet aisemment du th\'eor\`eme de Federer Fleming le r\'esultat suivant 
\begin{prop}
\label{masse-varifolds}
Soit $C_n$ une suite de courants int\'egraux d'une vari\'et\'e riemannienne compacte $(M^m,g)$ satisfaisant
\[
\limsup_{n\rightarrow +\infty} M(C_n)+M(\p C_n)<+\infty\quad.
\]
Alors il existe une sous-suite $C_{n'}$ qui converge faiblement vers un courant \underbar{int\'egral} $C_\infty$ et telle que $|C_{n'}|$ converge faiblement au sens
des mesures de Radon vers $V$. On a
 \[
M(C_{n'})\le \|V\|(M^m)=\lim_{n'\rightarrow+\infty}\||C_{n'}|\|(M^m)\quad,
\] et par ailleurs on a l'\'equivalence entre les 3 affirmations suivantes
\begin{itemize}
\item[i)]
\[
\lim_{n'\rightarrow +\infty} M(C_{n'})=M(C_\infty)\quad,
\]
\item[ii)]
\[
M(C_\infty)=\|V\|(M^m)\quad,
\]
\item[iii)]
\[
|C_\infty|=V\quad.
\]
\end{itemize}
\end{prop}

Avec pour objectif de mettre en \'evidence des points critiques de la masse pour des {\it varifolds}, il est naturel d'\'etudier les variations premi\`ere de celle ci.
Un diff\'eomorphisme $\phi$ de $M^m$ dans lui m\^eme d\'efini de fa\c on unique un diff\'eomorphisme $\Phi$ du fibr\'e en grassmanienne associ\'e $G_k(M^m)$
tel que tout $k-$plan $P_x$ en un point $x$ quelconque soit envoy\'e sur le $k-$plan des vecteurs $d\phi_xX$ o\`u $X\in P_x$.
Le {\it pouss\'e en avant} d'un varifold $V$ par $\phi$ est donc d\'efini naturellement par
\[
\forall A\ \mbox{ensemble borelien de }G_k(M^m)\quad (\phi_\ast V)(A):=\int_{\Phi^{-1}(A)} J_P\phi\ dV(P)
\]
o\`u pour tout $k-$plan non orient\'e $P$ de $T_xM^m$ on d\'efinit
\[
J_P\phi:=\sqrt{\det\lf[(d\phi_{|_P})^\ast (d\phi_{|_P})\rg]}=|d\phi(e_1)\wedge\cdots d\phi(e_k)|_g
\]
o\`u $(e_1\cdots e_k)$ est une base orthonorm\'ee de $P$ pour la m\'etrique $g$. Cette d\'efinition du {\it pouss\'e en avant} d'un {\it varifold} est naturel. En effet on v\'erifie aisemment
que si $C$ est un {\it courant entier rectifiable}, pour tout diffeomorphisme $\phi$ de $M^m$, si on note $\phi_\ast C$ son pouss\'e en avant d\'efini par
\[
\forall\om\in C^\infty_0(\wedge^{k}M^m)\quad\quad\left<\phi_\ast C,\om\right>:=\left<C,\phi^\ast\om\right>
\]
o\`u $\phi^\ast\om$ est la notation habituelle pour le {\it tir\'e en arri\`ere} par $\phi$ de la forme $\om$,  alors on a
\[
|\phi_\ast C|=\phi_\ast|C|
\]
Soit maintenant $X$ un champs de vecteur $C^1$ diff\'erentiable sur $M^m$. Nous notons ${\mathcal X}(M^m)$ l'espace de ces champs de vecteurs. Soit $\phi_t$ le flot sur $M^m$ associ\'e \`a ce champs de vecteur tel
que $\phi_0(x)=x$ pour tout $x\in M^m$. \'Etant donn\'e un
{\it varifold} $V$, la premi\`ere variation de $V$ par rapport \`a $X$ est d\'efinie par
\[
\delta V(X):=\lf.\frac{d}{dt}\rg|_{t=0}\|(\phi_t)_\ast V\|(M)
\]
Un {\it varifold} est dit {\it stationnaire} si
\[
\forall X\in {\mathcal X}(M)\quad\quad\delta V(X)=0\quad.
\]
Un r\'esultat important de W.Allard \cite{All} affirme que dans ${\R}^m$ un {\it varifold rectifiable entier stationnaire} est le {\it varifold} g\'en\'er\'e par une sous vari\'et\'e $C^\infty$ affect\'e d'une densit\'e enti\`ere r\'eguli\`ere au voisinage de tout point o\`u la densit\'e
est localement suffisamment proche d'une constante. On en d\'eduit qu'un tel {\it varifold} est r\'egulier dans un ouvert dense de son support. Le prob\`eme de d\'ecrire la r\'egularit\'e d'un {\it varifold rectifiable entier stationnaire} au voisinage de points o\`u la densit\'e varie est compl\`etement ouvert m\^eme dans le cas le plus simple dans {\it varifold rectifiable entier} bidimensionnel dans ${\R}^3$.

Dans la lign\'ee des r\'esultats de min-max sur les g\'eod\'esiques, pr\'esent\'es dans la premi\`ere partie de cet expos\'e, on va d\'efinir des classes d'homotopies d'applications continues \`a valeur dans les {\it cycles rectifiables entiers} ${\mathcal Z}_k(M^m)$. Pour cela il nous faut une topologie sur ces cycles. Nous avons d\'eja la topologie faible - celle des distributions - en dualit\'e avec les formes $C^\infty$ \`a support compacte mentionn\'ee plus haut en particulier dans le r\'esultat de fermeture des {\it courants int\'egraux} de Federer et Fleming . Cette topologie est en fait trop faible. Il se trouve que H.Federer et W.Fleming d\'emontrent une convergence plus forte de la sous suite $C_{n'}$ dans le th\'eor\`eme~\ref{FedererFleming} : la convergence pour la {\it topologie b\'emol}\footnote{L'origine de cette d\'enomination n'est pas connue de l'auteur de cet expos\'e. Cela \'etant il semble que l'image musicale pourrait se r\'ef\'erer au fait que la topologie
initialement consid\'er\'ee \'etait originellement celle de la {\it masse} - que nous pr\'esentons plus loin dans cet expos\'e - qui est une topologie plus forte ou plus ''haute'' disons et que donc la topologie dite
{\it b\'emol} est un affaiblissement  ou abaissement de celle ci....} . 

\medskip

La {\it topologie b\'emol} sur les {\it courants rectifiables entiers}\footnote{Elle s'\'etend en fait sur un espace plus grand qui est celui des
{\it courants rectifiables entiers}, l'espace des {\it courants b\'emol} : la compl\'etion des courants de masse finie et dont le bord a une masse finie pour la {\it distance b\'emol}, voir \cite{Fe}.} est une topologie d'espace vectoriel norm\'e issue
de la {\it norme b\'emol}, d\'efinie par
\[
\begin{array}{rcl}
{\mathcal F}(C):&=&\sup\lf\{<C,\om>\ ;\ \mbox{ t.q. } \|\om\|_\ast\le 1\ \mbox{ et }\ \|d\om\|_\ast\le 1\rg\}\\[5mm]
 &=&\inf\lf\{ M(A)+M(B)\ ;\ C= A+\p B\ \mbox{ t.q. } A\in{\mathcal D}_k(M)\ \mbox{ et }\ B\in{\mathcal D}_{k+1}(M)\rg\}
\end{array}
\]
La distance {\it b\'emol} entre deux {\it courants rectifiables entiers} $C_1$ et $C_2$ est donn\'ee par ${\mathcal F}(C_1,C_2):={\mathcal F}(C_1-C_2)$.
Par exemple si $C_1$ et $C_2$ sont des   {\it courants rectifiables entiers} de dimension 0, c'est \`a dire des sommes finies de masses de Dirac affect\'ees de multiplicit\'es
enti\`eres, et si on fait l'hypoth\`ese de neutralit\'e $C_i({\mathbf 1}_M)=0$, o\`u ${\mathbf 1}_M$ est la fonction identiquement \'egale \`a 1 sur $M$, dans le cas o\`u les supports
de $C_1$ et $C_2$ sont suffisamment proches, ${\mathcal F}(C_1,C_2)$ n'est rien d'autre que la {\it 1-distance de Waserstein} du transport optimal entre les deux familles de points affect\'es de multiplicit\'es
enti\`ere, c'est \`a dire la longueur minimale n\'ec\'essaire pour connecter ces deux familles de points affect\'es de ces multiplicit\'es.  Lorsque l'espace des {\it cycles rectifiable entier} de dimension $k$ est \'equip\'e de cette topologie on utilisera la notation ${\mathcal Z}_k(M,{\mathcal F})$.

\medskip

L'espoir de d\'evelopper une m\'ethode de {\it min-max} reprenant la m\'ethode de ''balayage'' de Birkhoff dans le cadre des {\it cycles rectifiable entier} est construit sur le r\'esultat suivant.
\begin{theo} [\cite{Alm}]
\label{Alm}
Pour tout $1\le n\le m$ et tout $1\le k\le m-n$  l'espace $H_{n+k}(M^m,{\Z})$ est isomorphe \`a $\pi_{n}({\mathcal Z}_k(M,{\mathcal F}))$.
\end{theo}

\medskip

Afin de g\'en\'erer des surfaces minimales non n\'ecessairement minimisantes, l'id\'ee serait alors de regarder un probl\`eme de {\it min-max}  de la forme
\be
\label{min-max-alm}
\inf_{\Phi\simeq\Psi}\sup_{x\in [0,1]^{n}}M(\Phi(x))
\ee
o\`u $\Psi$ r\'ealise une classe d'homotopie non triviale de ${\mathcal Z}_k(M,{\mathcal F})$ pour la donn\'ee au bord fix\'ee \'egale au courant nul\footnote{Cette condition de bord correspond g\'eom\'etriquement \`a l'hypoth\`ese de bord des balayages de Birkhoff. }et $\Phi$ \'evolue parmi toutes les applications continues
de $ [0,1]^{n}$ dans ${\mathcal Z}_k(M,{\mathcal F})$ qui sont homotopes \`a $\Psi$ pour des d\'eformations continues  coincidant avec $0$ sur le bord du cube $\p [0,1]^{n}$.

\medskip

Sh\'ematiquement il reste deux difficult\'es principales \`a surmonter avant de pouvoir r\'ealiser avec succ\'es un tel programe. 
\begin{itemize}
\item[i)] La masse n'est pas continue mais seulement semi-continue inf\'erieurement pour la {\it distance b\'emol} et donc la topologie consid\'er\'ee pour
d\'efinir les classes d'homotopie ne devrait pas permettre \`a priori d'obtenir la r\'ealisation de l'infimum, c'est \`a dire, avec la topologie {\it b\'emol} il peut tr\`es bien y  avoir un ''saut de masse'' par
passage \`a la limite dans (\ref{min-max-alm}).
\item[ii)] Si jamais un tel {\it min-max} (\ref{min-max-alm}) \'etait atteint et r\'ealis\'e par un $\Phi(x_0)$, qu'en est-il de sa r\'egularit\'e ? Nous aurions peut \^etre bien g\'en\'er\'e ainsi un {\it varifold rectifiable entier stationaire} en consid\'erant  $|\Phi(x_0)|$ mais le r\'esultat de W.Allard n'est pas assez fort pour en d\'eduire la r\'egularit\'e esp\'er\'ee et pour affirmer que l'on a bien affaire \`a une immersion minimale r\'eguli\`ere. 
\end{itemize}

\medskip

L'\'etudiant de F.Almgren, J.Pitts, s'est attel\'e \`a la r\'esolution de ces deux principales difficult\'es et nous pr\'esentons ci-dessous les grandes lignes de son travail spectaculaire, \cite{Pit},
qui lui a permis en particulier de g\'en\'erer des familles de surfaces minimales jusque l\`a inconnues.

Le projet de F.Almgren, d\'eja bien conscient des difficult\'es i) et ii) mentionn\'ees plus haut,  puis de Pitts, \'etait de d\'evelopper une th\'eorie de {\it min-max} pour les {\it cycles rectifiables entiers} ''accompagn\'es'' de leur {\it varifold} associ\'e afin de ne rien perdre
de la masse et que l'''objet'' obtenu - on \'esp\`ere un {\it varifold  rectifiable entier stable} -  soit bien ''accroch\'e'' au sens o\`u sa masse r\'ealise la {\it largeur} du {\it min-max} .
 J. Pitts \cite{Pit} a  donc propos\'e de ''renforcer'' la {\it topologie b\'emol} pour rem\'edier \`a cette difficult\'e. Sur l'espace des {\it varifolds rectifiables entiers} tout d'abord on d\'efinit la m\'etrique suivante
\[
{\mathbf F}(V,W):=\sup\lf\{V(f)-W(f)\ ;\ f\in C_c(M^m)\ \mbox{ t. q. }|f|\le 1\ \mbox{ et } Lip(f)\le 1\rg\}
\] 
o\`u $Lip(f)=\sup_{x\ne y}|f(x)-f(y)|/d_M(x,y)$ est la norme lipschitz de $f$. Sur l'espace des {\it courants rectifiables entiers} on d\'efinit alors la ${\mathbf F}-$m\'etrique
\[
{\mathbf F}(C_1,C_2):={\mathcal F}(C_1-C_2)+{\mathbf F}(|C_1|,|C_2|)\quad.
\]
On v\'erifie aisemment que la masse est continue pour la topologie induite par la distance ${\mathbf F}$.
\medskip

Malheureusement nous ne sommes pas au bout de nos peines pour avoir la bonne topologie pour laquelle le probl\`eme de {\it min-max} sera bien pos\'e. La topologie la plus forte 
dans le cadre des {\it courants rectifiables entiers} est la topologie dite {\it topologie de la masse} d\'efinie par la distance suivante
\[
\forall \ C_1, C_2\in {\mathcal R}_k(M)\quad{\mathbf M}(C_1,C_2):=M(C_1-C_2)
\]
Il n'est pas difficile de v\'erifier l'emboitement suivant des 3 topologies ci dessus.
\[
{\mathcal F}(C_1,C_2)\le {\mathbf F}(C_1,C_2)\le 2\, {\mathbf M}(C_1,C_2)\quad.
\]
Lorsque l'espace des {\it cycles rectifiable entier} de dimension $k$ est \'equip\'e de chacune de ces 3 topologies on utilisera les notations 
respectivement ${\mathcal Z}_k(M,{\mathcal F})$, ${\mathcal Z}_k(M,{\mathbf F})$ et ${\mathcal Z}_k(M,{\mathbf M})$.

\medskip

Inspir\'e par le travail de M.Morse \cite{Mo} sur les m\'ethodes topologiques dans le calcul des variations, plut\^ot que de consid\'erer des homotopies continues d'applications de
$[0,1]^n$ dans ${\mathcal Z}_k(M,{\mathbf F})$, J.Pitts va lui pr\'ef\'erer des versions discr\`etes \`a valeur dans ${\mathcal Z}_k(M,{\mathbf M})$ qui sont \`a la fois plus ''informatives'' et plus souples pour proc\'eder \`a des arguments de comparaisons cela
afin de d\'emontrer la criticalit\'e et l'{\it approximative minimalit\'e}, que nous d\'efinissons plus bas, du maximum de l'application limite.

\medskip

L'introduction des suites d'homotopie discr\`etes requiert un peu de notation. Ces notions sont tr\`es naturelles mais leur pr\'esentation n'\'echappe malheureusement pas \`a une certaine lourdeur.

On note $I=[0,1]$ et pour tout $j\in {\N}$ on note $I(1,j)$ le complexe cellulaire dont les 
cellules de dimensions $0$ sont les points $[k\,3^{-j}]$ pour $k=0,1\cdots 3^j$ et les cellules de dimension 1 sont les segments $[k\,3^{-j},(k+1)\,3^{-j}]$. Plus g\'en\'eralement
$I(n,j)$ est le complexe cellulaire dont les cellules $p-$dimensionnelles s'identifient aux sous ensembles de $[0,1]^n$ de la forme $\al_1\times\cdots\times\al_n$ telles
que $\sum_{i=1}^n\dim \al_i=p$ o\`u $\al_i$ sont soit des cellules de dimension $0$ soit de dimension $1$ dans $I(1,j)$. L'espace des sommets de ce complexe, qui sera not\'e $I(n,j)_0$, s'identifie 
donc \`a l'espace des points de coordonn\'ees $(k_1\, 3^{-j},\cdots, k_n\, 3^{-j})$ pour un choix quelconque d'entiers $k_1\cdots k_n$ parmi $\{0,\cdots, 3^j\}$.

Par un abus de notation on \'ecrit $\p I(n,j)_0$ l'intersection entre $I(n,j)_0$ et le bord du cube $\p [0,1]^n$.

La distance
entre deux sommets $x=(k_1\, 3^{-j},\cdots, k_n\, 3^{-j})$ et $x'=(k'_1\, 3^{-j},\cdots, k'_n\, 3^{-j})$ est donn\'ee par
\[
{\mathbf d}_j(x,x'):=\sum_{i=1}^n|k_i-k'_i|
\]
\'Etant donn\'ee une application $\varphi^j$ de $I(n,j)_0$ dans ${\mathcal Z}_k(M)$, on d\'efinie la {\it finesse} de $\varphi^j$ par 
\[
{\mathbf f}(\varphi^j):=\sup\lf\{ \frac{{\mathbf M}(\varphi^j(x)-\varphi^j(x'))}{{\mathbf d}_j(x,x')}\ ;\ \forall\ x\, ,\,x'\ \in I(n,j)_0\ \mbox{ et } x\ne x'\rg\}
\]

\medskip

On cherche maintenant \`a d\'efinir une notion ''d'homotopie discr\`ete'' entre deux applications $\varphi^{j_1}_1$ et $\varphi^{j_2}_2$ de $I(n,j)_0$ dans ${\mathcal Z}_k(M)$
toutes deux \'egales \`a zero sur le bord $\p I(n,j)_0$.

 Lorsque $j_1=j_2=j$ tout d'abord,
pour tout $\delta>0$, on dite que $\varphi_1^j$ est {\it homotope} \`a $\varphi_2^j$ pour la {\it finesse} $\delta$ si il existe
\[
\psi^j\ :\ I(1,j)_0\times I(n,j)_0\ \longrightarrow {\mathcal Z}_k(M)
\]
telle que
\begin{itemize}
\item[i)]
\[
\forall \, x\in I(n,j)_0\ \quad\quad\psi^j(0,x)=\varphi^{j_1}_1(x)\quad\mbox{ et }\quad \psi^j(1,x)=\varphi^{j_2}_2(x)
\]
\item[ii)]
\[
\forall\, y\in I(1,j)_0\times \p I(n,j)_0\quad \quad\quad\psi^j(y)=0
\]
\item[iii)]
\[
{\mathbf f}(\psi^j)<\delta
\]
\end{itemize}
Lorsque maintenant $j_1<j_2$, on g\'en\`ere alors une application $\varphi_1^{j_2}$ \`a partir de $\varphi_1^{j_2}$ de la fa\c con suivante :
pour tout $x\in I(n,j_2)_0$ on prend 
\[
\varphi_1^{j_2}(x):=\varphi_1^{j_1}(n(j_2,j_1)(x))
\] o\`u $n(j_2,j_1)(x)$ est le sommet du complexe $I(n,j_1)_0$, inclus
dans le complexe $I(n,j_2)_0$ qui minimize la distance ${\mathbf d}_{j_2}$. 

On dit alors que $\varphi_1^{j_1}$ est {\it homotope} \`a $\varphi_2^{j_2}$ pour la {\it finesse} $\delta$ si
$\varphi_1^{j_2}$ est {\it homotope} \`a $\varphi_2^{j_2}$ pour la {\it finesse} $\delta$.

On introduit alors les deux d\'efinitions suivantes qui sont au coeur de la construction de Pitts.
\begin{defi}
\label{n-m-suite-homotopie}
Soit $n\in {\N}^\ast$. Une suite $\varphi=\{\varphi^j\}_{j\in{N}}$ o\`u chaque $\varphi^j$ est une application de $I(n,j)_0$ dans ${\mathcal Z}_k(M)$ est appel\'ee 
une $(n,{\mathbf M})$ - suite homotopique \`a valeur dans ${\mathcal Z}_k(M)$ si il existe une suite de r\'eels strictement positifs $\delta_i$
telle que
\begin{itemize}
\item[i)]
pour tout $j\in {N}$ $\varphi^j$ est homotope \`a $\varphi^{j+1}$ pour la finesse $\delta^j$.
\item[ii)]
\[
\lim_{j\rightarrow +\infty}\delta^j=0
\]
\item[iii)]
\[
\sup_{j\in{\N}}\ \sup_{x\in I(n,j)_0}{\mathbf M}(\varphi^j(x))<+\infty\quad.
\]

\end{itemize}
\end{defi}

La deuxi\`eme d\'efinition introduit la notion de classe d'homotopie discr\`ete.

\begin{defi}
\label{n-m-classe-homotopie}
Soient $\varphi_1=\{\varphi_1^j\}_{j\in{N}}$ et $\varphi_2=\{\varphi_2^j\}_{j\in{N}}$  deux $(n,{\mathbf M})$ - suites homotopiques \`a valeur dans ${\mathcal Z}_k(M)$ sont
homotopes si il existe une suite de r\'eels positifs $\delta^j>0$ telle que
\begin{itemize}
\item[i)]
 $\varphi_1^j$ est homotope \`a $\varphi^{j}_2$ pour la finesse $\delta^j$ pour tout $j\in {\N}$.
\item[ii)]
\[
\lim_{j\rightarrow +\infty}\delta^j=0
\]
\end{itemize}
C'est une relation d'\'equivalence dont l'espace des classes est not\'ee $\pi_n^\sharp({\mathcal Z}_k(M, {\mathbf M}))$. La classe d'une $(n,{\mathbf M})$ - suites homotopique $\varphi$
sera not\'ee $[\varphi]$.
\end{defi}
Le r\'esultat suivant est une extension du r\'esultat d'Algmgren pr\'esent\'e plus haut.

\begin{theo} [\cite{Pit}]
\label{Alm-Pitt}
Pour tout $1\le n\le m$ et tout $1\le k\le m-n$  l'espace $H_{n+k}(M^m,{\Z})$ est isomorphe \`a $\pi_{n}({\mathcal Z}_k(M,{\mathcal F}))$ et \`a $\pi_n^\sharp({\mathcal Z}_k(M, {\mathbf M}))$.
\end{theo}
Ce r\'esultat peut para\^ \i tre au premier abord tout \`a fait surprenant. En effet la continuit\'e pour la topologie $\mathcal F$ semble beaucoup plus faible que la continuit\'e discr\`ete
avec une finesse tendant vers z\'ero pour la topologie plus forte ${\mathbf M}$. On s'attendrait donc \`a avoir beaucoup plus de classes d'homotopie
pour la deuxi\`eme topologie. La figure suivante devrait pouvoir \'eclairer ce paradoxe et rendre le th\'eor\`eme pr\'ec\'edent plus naturel.

\psfrag{a}{\'etape 1}
\psfrag{b}{\'etape 2}
\psfrag{c}{\'etape 3}
\psfrag{aa}{\'etape $p^2 -2$}
\psfrag{bb}{\'etape $p^2-1$}
\psfrag{cc}{\'etape $p^2$}
\psfrag{f}{finesse}
\psfrag{d}{$\delta  = 1/p$}
\psfrag{h}{Homotopie ${\mathcal F}$-continue}
\psfrag{p}{passsant de ~~~~~ \`a}
\psfrag{A}{$(1,{\mathbf{M}})$-homotopie}
\psfrag{q}{passsant de ~~~~~ \`a}
\psfrag{g}{pour la finesse $\delta$}
\begin{center}
\includegraphics[width=8cm]{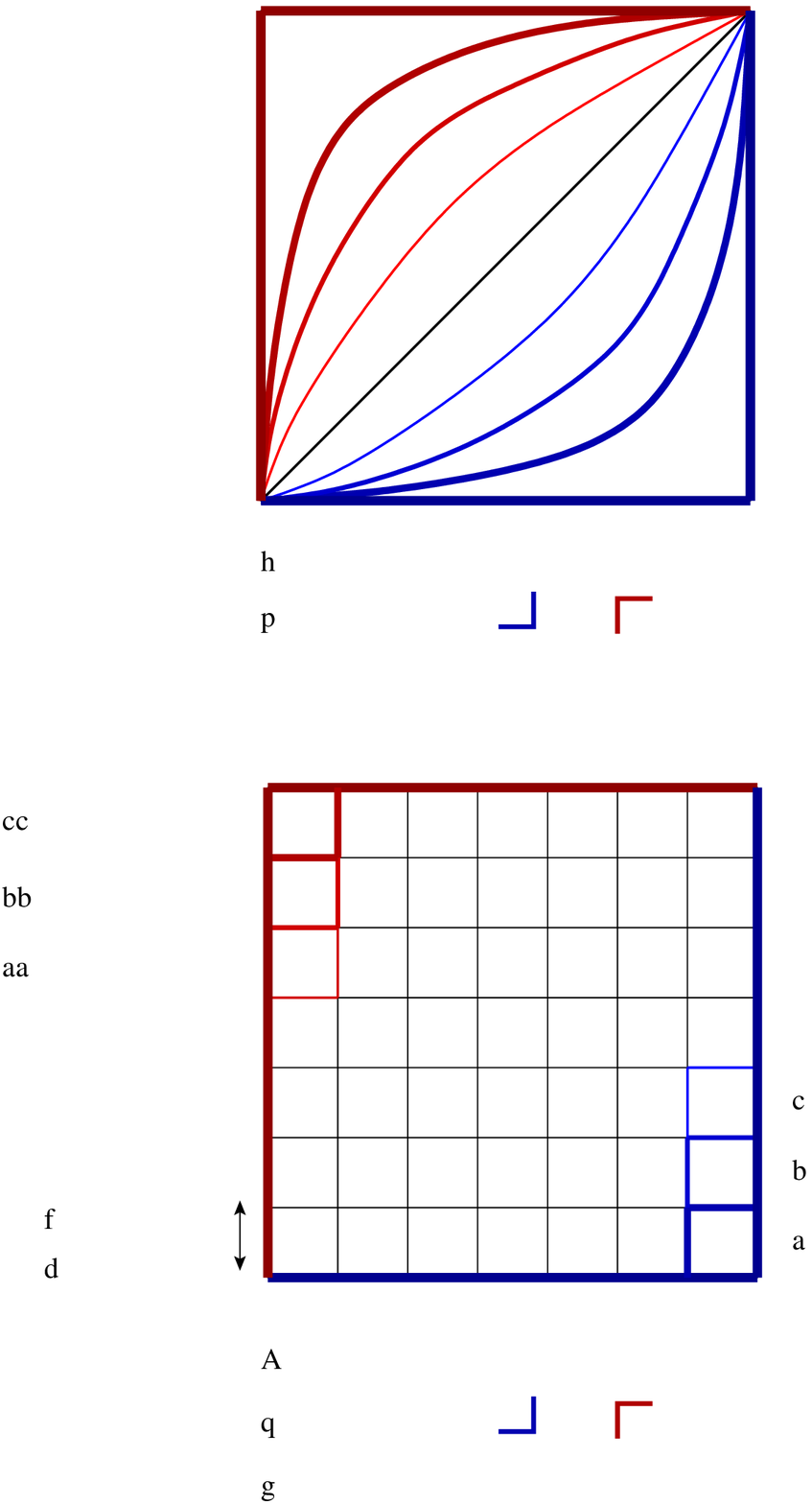}
\end{center}

\medskip

On se propose d'\'etudier des probl\`emes de {\it min-max} \`a partir de balayage de $M$ r\'ealis\'es par des classes non nulles de $\pi_n^\sharp({\mathcal Z}_k(M, {\mathbf M}))$. Soit donc $\varphi=(\varphi^j)$ une $(n,{\mathbf M})$ - {\it suite homotopique }\`a valeur dans ${\mathcal Z}_k(M)$ on d\'efinit
\[
{\mathbf L}(\varphi):=\limsup_{j\rightarrow+\infty}\max\lf\{ {\mathbf M}(\varphi^j(x))\ ;\ x\in I(n,j)_0\rg\}\quad.
\]
C'est le remplacement en version discr\`ete du maximum de la masse d'une application de $[0,1]^n$ dans ${\mathcal Z}_k(M)$. On d\'efinit alors la {\it largeur} du probl\`eme de {\it min-max} associ\'e \`a une classe de $\pi_n^\sharp({\mathcal Z}_k(M, {\mathbf M}))$ donn\'ee dans la continuit\'e des d\'efinitions de la partie 1 de cet expos\'e. 
\begin{defi}
\label{largeur}
Soit $n\in {\N}^\ast$ et $k\in {\N}^\ast$ et soit $\Pi\in \pi_n^\sharp({\mathcal Z}_k(M, {\mathbf M}))$. La largeur du min-max associ\'e \`a $\Pi$ est donn\'ee par
\[
{\mathbf L}(\Pi):=\inf\{{\mathbf L}(\varphi)\ ;\ [\varphi]\in \Pi\}
\]
On dit que $\varphi$ est critique pour $\Pi$ si ${\mathbf L}(\varphi)={\mathbf L}(\Pi)$.
\end{defi}
Soit $\varphi=(\varphi^j)$ une $(n,{\mathbf M})$ - {\it suite homotopique} \`a valeur dans ${\mathcal Z}_k(M)$ on d\'efinit l'espace limite de $\varphi$ dans ${\mathcal V}_k(M)$ par
\[
{\mathbf K}(\varphi):=\lf\{V= \lim_{j'\rightarrow+\infty}|\varphi^{j'}(x^{j'})| \mbox{ pour une sous-suite }j'\mbox{ et pour }x^{j'}\in I(n,j')_0\rg\}
\]
la limite ci-dessus \'etant au sens des mesures de Radon. Enfin on d\'efinit l'ensemble critique de $\varphi$ comme \'etant
\[
{\mathbf C}(\varphi):={\mathbf K}(\varphi)\cap \{V\ ;\ \|V\|(M)={\mathbf L}(\varphi)\}
\]
Il est important de noter \`a ce stade que $C(\varphi)$ est compacte et non-vide, c'est l\`a tout l'avantage de ''suivre'' la projection de  nos {\it cycles entiers rectifiables} dans l'espace des {\it varifolds}. C'est cette projection qui va contenir toute l'information de criticalit\'e et de la quasi-minimalit\'e - une version faible d'\^etre d'indice fini - de notre suite critique $\varphi$.

\medskip

Le premier succ\'es de l'approche de Pitts est contenu dans ce th\'eor\`eme.
\begin{theo} [\cite{Pit}]
\label{Pitt-exist-critical}
Soit $M^m$ une vari\'et\'e ferm\'ee. Soit $n\in {\N}^\ast$ et $k\in {\N}^\ast$. Soit $\Pi\in\pi_n^\sharp({\mathcal Z}_k(M, {\mathbf M}))$,  alors il existe une $(n,{\mathbf M})$ - {\it suite homotopique }\`a valeur dans ${\mathcal Z}_k(M)$, $\varphi_\ast$,
qui soit critique, c'est \`a dire telle que ${\mathbf L}(\varphi_\ast)={\mathbf L}(\Pi)$. Par ailleurs $\varphi_\ast$ peut \^etre choisie telle que chaque \'el\'ement de ${\mathbf C}(\varphi_\ast)\ne\emptyset$ soit un varifold stationnaire.
\end{theo}

Cela \'etant notre enthousiasme doit \^etre encore contenu \`a ce stade car il se pourrait que $L(\Pi)=0$ et que les \'el\'ements de $C(\varphi_\ast)$ soient triviaux. Le th\'eor\`eme suivant 
dont la preuve se d\'eduit de celle du th\'eor\`eme~\ref{Alm-Pitt} plus haut nous assure que ce n'est pas toujours le cas.

\begin{theo} [\cite{Pit}]
\label{Pitt-exist-critical-1}
Soit $M^m$ une vari\'et\'e ferm\'ee et soit $1\le n\le m$ et  $1\le k\le m-n$  tels que $H_{n+k}(M^m,{\Z})\ne 0$ alors il existe une classe $\Pi\in \pi_n^\sharp({\mathcal Z}_k(M, {\mathbf M}))$ telle que $L(\Pi)>0$.
\end{theo}

Si donc on combine les deux th\'eor\`emes pr\'ec\'edents on arrive \`a g\'en\'erer un {\it varifold stationnaire} non-trivial d\`es lors que l'on a un groupe $H_{n+k}(M^m,{\Z})$ non trivial. Que faire alors d'un tel objet ? le r\'esultat d'Allard ne nous assure la r\'egularit\'e de ce {\it varifold} sur un ensemble ouvert dense, et cela  seulement si  on sait que ce {\it varifold} est rectifiable, ce que l'on ne sait pas encore.  C'est sur ces questions de r\'egularit\'e que le travail de Pitts apporte une contribution fondamentale. L'id\'ee est que si une surface minimale r\'ealiserait un tel {\it min-max} alors elle serait d'indice fini et elle aurait la propri\'et\'e d'\^etre ''presque'' d'aire minimale globalement
et de l'\^etre certainement localement. Or nous avons vu qu'il existe des r\'esultats de r\'egularit\'e tr\`es fort pour les cycles rectifiables d'aire minimale et de codimension 1. 
Pitts va donc formaliser cette notion de ''presque minimalit\'e'' pour les {\it varifolds} et plus de la moiti\'e de son travail sera d'en d\'emontrer la r\'egularit\'e en codimension 1 pour toute dimension de l'espace ambiant allant de 3 \`a 7.

Pour d\'efinir les {\it varifolds presques minimaux} il nous faut \'etendre quelques notations. \'Etant donn\'e un sous ensemble ferm\'e $N$ de la vari\'et\'e $M^m$ on note ${\mathcal Z}_k(M,N)$ l'espace des courants rectifiables entiers dont le bord est support\'e dans $N$. \'Etant donn\'e un ouvert $U$ de $M$ on note ${\mathbf F}_U$ la semi distance sur les {\it varifolds} obtenue en restreignant la mesure de la distance ${\mathbf F}$ aux parties de ces {\it varifolds} contenues dans $U$. On a alors la d\'efinition suivante.
\begin{defi}
\label{presque-minim}
Un varifold $V$ dans ${\mathcal V}_k(M^m)$ et dit presque minimisant dans un ouvert $U$ de $M^m$ si pour tout $\epsilon>0$ il existe $\delta>0$ tels que
pour tout \'el\'ement $C$ de ${\mathcal Z}_k(M,M\setminus U)$ satisfaisant ${\mathbf F}_U(V,|C|)<\ep$ et toute suite finie $(C_i)_{i=1\cdots p}$ satisfaisant $C_1=C$ et
\[
\lf\{
\begin{array}{l}
\ds supp(C-C_i)\subset U\quad\quad\ \ \forall i=2\cdots p\\[5mm]
\ds {\mathbf M}(C_{i}-C_{i-1})\le \delta\ \ \quad\quad\forall i=2\cdots p\\[5mm]
\ds {\mathbf M}(C_i)\le M(C)+\delta\quad\quad\forall i=2\cdots p
\end{array}
\rg.
\] 
alors
\[
{\mathbf M}(C_p)\ge M(C)-\ep
\]
\end{defi}
Cette d\'efinition peut se comprendre comme suit. Un {\it varifold} $V$ est {\it presque minimisant} si, pour tout $\epsilon$, si on cherche \`a ''connecter'' un courant  $C$  \`a  ${\mathbf F}$ distance $\epsilon$ de $V$ et un courant $C_p$ de masse strictement plus petite que ${\mathbf M}(C)-\epsilon$ par une cha\^ \i ne finie de courants $C_i$ telle que $M(C_i- C_{i-1})\le\delta$ alors il existe n\'ec\'essairement
un de ceux-ci qui a une masse sup\'erieure d'au moins $\delta$ de ${\mathbf M}(C)$. En termes simplifi\'es cela implique que toute famille discr\`ete $\delta-${\it fine}, coincidant avec $V$ en dehors de $U$ et dont un \'el\'ement  a une masse
plus basse que $\|V\|(U)-\epsilon$ alors un de ses \'el\'ements \`a une masse plus \'elev\'ee que $\|V\|(U)+\delta$ et donc le {\it max} de cette famille est plus \'elev\'e que $\|V\|(M)$. On comprend donc bien comment cette propri\'et\'e sera satisfaite par les \'el\'ements 
critiques de $C(\varphi_\ast)$ r\'ealisant le {\it min-max}. Le r\'esultat suivant vient en partie le confirmer.

\begin{theo} [\cite{Pit}]
\label{Pitt-exist-critical-stationnaire}
Soit $M^m$ une vari\'et\'e ferm\'ee. Soit $n\in {\N}^\ast$ et $k\in {\N}^\ast$. Soit $\Pi\in\pi_n^\sharp({\mathcal Z}_k(M, {\mathbf M}))$,  alors il existe une $(n,{\mathbf M})$ - {\it suite homotopique }\`a valeur dans ${\mathcal Z}_k(M)$, $\varphi_\ast$,
qui soit critique, c'est \`a dire telle que ${\mathbf L}(\varphi_\ast)={\mathbf L}(\Pi)$. Par ailleurs il existe un \'el\'ement de ${\mathbf C}(\varphi_\ast)$ qui soit un varifold stationnaire
et presque minimisant sur tout anneau ouvert.\end{theo}

Vient ensuite un travail de r\'egularit\'e sur les {\it varifolds stationnaires} de ${\mathcal V}_k(M)$  et on d\'emontre dans un premier temps  le r\'esultat suivant.

\begin{theo} [\cite{Pit}]
Tout \'el\'ement de ${\mathcal V}_k(M)$ qui soit {\it stationnaire} et {\it presque minimisant} sur tout anneau ouvert est dans 
l'espace des {\it varifolds int\'egraux} $IV_k(M)$.
\end{theo}

Le th\'eor\`eme de r\'egularit\'e pour les {\it varifolds presque minimisant} de \underbar{codimension 1} est d\'emontr\'e tout d'abord dans \cite{Pit} pour $m\le 6$ puis dans \cite{ScS} pour $m=7$. Dans ce dernier travail R.Schoen et L.Simon \'etendent ce r\'esultat en un r\'esultat de r\'egularit\'e partielle en dimension quelconque. La discussion plus haut explique pourquoi \`a partir de la dimension 8 il faut s'attendre \`a ce que de tels  {\it varifolds} aient des singularit\'es.

\begin{theo} 
Soit $M^m$ une vari\'et\'e ferm\'ee de dimension $m\le 7$.  Tout {\it varifold int\'egral stationnaire}  de codimension 1 et {\it presque minimisant} sur tout anneau ouvert est le varifold associ\'e 
au courant d'int\'egration le long d'une sous vari\'et\'e minimale r\'eguli\`ere affect\'ee d'une multiplicit\'e enti\`ere r\'eguli\`ere.
\end{theo}

En combinant le travail de {\it min-max} plus haut et ce dernier r\'esultat de r\'egularit\'e sur les  {\it varifolds presque minimisants} de {codimension 1} on obtient le r\'esultat qui est le point
culminant du travail de Pitts.
\begin{theo} [\cite{Pit}]
\label{Pitt-exist-critical-codim-1}
Soit $M^m$ une vari\'et\'e ferm\'ee de dimension $m\le 7$. Soit $n\in {\N}^\ast$ et Soit $\Pi\in\pi_n^\sharp({\mathcal Z}_{m-1}(M, {\mathbf M}))$,  alors il existe une $(n,{\mathbf M})$ - {\it suite homotopique }\`a valeur dans ${\mathcal Z}_k(M)$, $\varphi_\ast$,
qui soit critique, c'est \`a dire telle que ${\mathbf L}(\varphi_\ast)={\mathbf L}(\Pi)$. Par ailleurs, si ${\mathbf L}(\Pi)>0$, il existe un \'el\'ement de ${\mathbf C}(\varphi_\ast)$ qui soit le varifold associ\'e 
au courant d'int\'egration le long d'une sous vari\'et\'e minimale r\'eguli\`ere de codimension 1 affect\'ee d'une multiplicit\'e enti\`ere r\'eguli\`ere.
\end{theo}
En combinant th\'eor\`eme~\ref{Pitt-exist-critical-1} et le r\'esultat pr\'ec\'edent on obtient finalement le corollaire suivant.

\begin{coro}[\cite{Pit}]
\label{exist-minimal}
Toute vari\'et\'e ferm\'ee r\'eguli\`ere de dimension inf\'erieure ou \'egale \`a $7$ poss\'ede une hyper-surface minimale plong\'ee.
\end{coro}

\section{La d\'emonstration de la conjecture de Willmore par F.Marques et A.Neves.}

\subsection{La conjecture de Willmore.}

Soit $\Sigma$ une vari\'et\'e bi-dimensionnelle orient\'ee et $\vec{\Phi}$ une immersion de cette vari\'et\'e dans une vari\'et\'e riemannienne $(M^m,g)$. Cette immersion induit une m\'etrique sur 
$\Sigma$ - la {\it premi\`ere forme fondamentale} de l'immersion - que nous notons $g_{\vec{\Phi}}$ qui est \'egale au {\it tir\'e en arri\`ere} par $\vec{\Phi}$ de la m\'etrique $g$ : 
$g_{\vec{\Phi}}(X,Y):=g(\vec{\Phi}_\ast X,\vec{\Phi}_\ast Y)$ o\`u $\vec{\Phi}_\ast X$ et $\vec{\Phi}_\ast Y$ sont les {\it pouss\'es en avant} de deux vecteur arbitraires $X$ et $Y$ tangents en un point quelconque de $\Sigma$. Parfois, lorsqu'il n'y a pas d'ambiguit\'e sur l'immersion consid\'er\'ee ces {\it pouss\'es en avant} seront simplement not\'es $\vec{X}$ et $\vec{Y}$.

Cette m\'etrique induit une forme volume sur $\Sigma$ que nous notons $dvol_{g_{\vec{\Phi}}}$ et qui est donn\'ee en coordonn\'ees locales positives $x=(x_1,x_2)$ par
$dvol_{g_{\vec{\Phi}}}= \sqrt{g_{11}g_{22}-g_{12}^2}\ dx_1\wedge dx_2$  o\`u  $g_{\vec{\Phi}}= \sum_{i,j=1}^2 g_{ij}\ dx_i\otimes dx_j\ $.

La premi\`ere variation de la fonctionnelle d'aire, \'egale \`a l'int\'egrale de la forme volume sur $\Sigma$ : Aire$(\vec{\Phi}):=\int_{\Sigma}dvol_{g_{\vec{\Phi}}}$, est donn\'ee par
\be
\label{vol}
\lf.\frac{d \mbox{Aire}(\vec{\Phi}+t\,\vec{w})}{dt}\rg|_{t=0}=-\,2\,\int_\Sigma\vec{H}_{\vec{\Phi}}\cdot \vec{w}\ dvol_{g_{\vec{\Phi}}}\quad,
\ee
o\`u $\vec{H}_{\vec{\Phi}}$ est le {\it vecteur courbure moyenne} de l'immersion $\vec{\Phi}$ qui est \'egale \`a la moiti\'e de la trace par rapport \`a la m\'etrique $g_{\vec{\Phi}}$ de la {\it seconde forme fondamentale} $\vec{\mathbb I}_{\vec{\Phi}}$ de l'immersion : en coordonn\'ees locale on a $$\vec{H}_{\vec{\Phi}}:= \frac{1}{2} tr_g\vec{\mathbb I} =\frac{1}{2} \sum_{i,j=1}^2g^{ij}\ \vec{\mathbb I}_{\vec{\Phi}}(\p_{x_i},\p_{x_j})$$ o\`u $(g^{ij})_{ij=1,2}$ est l'inverse de la matrice $(g_{ij})_{ij=1,2}$. On rappelle que la seconde forme fondamentale de $\vec{\Phi}$ est un \underbar{tenseur} qui \`a une paire de vecteurs $X$ et $Y$ tangent en un point $p$  quelconque de  $\Sigma$
associe un vecteur de $T_{\vec{\Phi}(p)}M^m$ orthogonal \`a la surface $\vec{\Phi}(\Sigma)$ et donn\'e par
\[
\vec{\mathbb I}(X,Y):=\pi_{\vec{n}}\lf(\nabla^g_Y(\vec{\Phi}_\ast X)\rg)=\pi_{\vec{n}}\lf(\nabla^g_X(\vec{\Phi}_\ast Y)\rg)=\vec{\mathbb I}(Y,X)
\]
o\`u $X$ et $Y$ sont \'etendus de fa\c con r\'eguli\`ere et arbitrare au voisinage de $p$, $\pi_{\vec{n}}$ est la projection orthogonale de $T_{\vec{\Phi}(p)}M$ dans le sous espace $[\vec{\Phi}_\ast(T_p\Sigma)]^\perp$ des vecteurs ortogonaux \`a $\vec{\Phi}_\ast(T_p\Sigma)$, enfin $\nabla^g$ est la {\it connection de Levi-Civita} de $(M^m,g)$ qui est de {\it torsion} nulle ce qui justifie la derni\`ere in\'egalit\'e et donc le fait que $\vec{\mathbb I}$ soit un tenseur symm\'etrique.

Comme nous l'avons vu plus haut les surfaces dites {\it minimales} sont les points critiques du volume et sont les g\'en\'eralisation naturelles des g\'eod\'esiques en dimension deux. L'expression (\ref{vol}), qui se g\'en\'eralise en dimension quelconque, nous dit qu'une immersion est {\it minimale} si et seulement si
\[
\vec{H}_{\vec{\Phi}}\equiv 0 \quad\quad\mbox{ sur }\Sigma\quad.
\]
L'\'etude des {\it surfaces minimales} est un domaine important des math\'ematiques qui peut s'aborder de fa\c cons tr\'es vari\'ees. L'approche variationnelle sur laqu'elle  nous avons insist\'e dans la premi\`ere partie de l'expos\'e, avec le probl\`eme
de Plateau et les m\'ethodes de {\it min-max}, est une des m\'ethodes d'\'etude de ces objets parmi de nombreuses autres comme les approches alg\'ebriques et g\'eom\'etriques - repr\'esentations de Weerstrass, repr\'esentations twistorielles...etc - ou les approches analytiques - Equations au d\'eriv\'ees partielles, m\'ethodes d'analyse fonctionnelle...etc. L'abondance des points de vues \'etant dus manifestement \`a l'universalit\'e de ces objets qui ressurgissent dans de nombreuses questions des math\'ematiques.

Dans un \'effort de ''fusioner'' la th\'eorie des surfaces minimales et l'invariance conforme - l'invariance par les transformations qui pr\'eservent les angles infinit\'esimalement - W.Blaschke au d\'ebut du XXeme si\`ecle \cite{Bl} \'etend l'espace des surfaces minimales \`a l'espace des points critiques de la fonctionnelle suivante
\[
W(\vec{\Phi}):=\int_\Sigma|\vec{H}_{\vec{\Phi}}|^2\ dvol_{g_{\vec{\Phi}}}\quad.
\]
Les surfaces minimales, qui sont de fa\c con \'evidente des minima absolus de ce {\it lagrangien }, n'en sont pas les seuls point critiques. On v\'erifie par exemple que la sph\`ere $S^2$ dans l'espace euclidien ${\R}^3$ est aussi point critique de $W$. 

W.Blaschke, qui travaillait dans le cadre de la g\'eom\'etrie conforme, proposa ce lagrangien car il observa une propri\'et\'e importante de celui-ci dans ${\R}^3$ : pour toute immersion $\vec{\Phi}$ d'une surface compacte ferm\'ee $\Sigma$- c'est \`a dire $\Sigma$ est compacte sans bord - et pour toute transformation conforme $\Psi$ de ${\R}^3\cup\{\infty\}$ dans ${\R}^3\cup\{\infty\}$ qui n'est singuli\`ere en aucun point de la surface - c'est \`a dire telle que  $\Psi^{-1}(\{\infty\})\cap \vec{\Phi}(\Sigma)=\emptyset$ - on a
\be
\label{invariance}
W(\Psi\circ\vec{\Phi})=W(\vec{\Phi})\quad.
\ee
on obtient ainsi en particulier que toute transformation conforme d'une {\it surface minimale} n'est peut-\^etre plus n\'ec\'essairement {\it minimale} mais est encore un point critique de $W$. C'est pourquoi probablement W.Blaschke d\'ecida de nommer ces surfaces les {\it surfaces minimales conformes}. L'origine de cette fonctionnelle remonte en fait \`a bien avant la g\'eom\'etrie conforme et le travail de Blaschke. Un si\`ecle plus t\^ot d\'ej\`a, dans son \'effort de g\'en\'eraliser la th\'eorie des poutres de J.Bernouilli et L.Euler aux membranes \'elastiques, S.Germain produit
un {\it lagrangien} qui fait intervenir la norme $L^2$ de la valeur moyenne en chaque point de toutes les courbures g\'eod\'esiques de l'intersection de la surface avec des plans perpendiculaire en ce point, qui n'est rien d'autre que $W$. Ses travaux\footnote{Ce qui est d'autant plus remarquable pour l'\'epoque car la g\'eom\'etrie diff\'erentielle des surfaces vers 1810 n'\'etait qu'a ses balbutiements et le {\it theorema egregium} de Gauss par exemple ne sera d\'ecouvert que 18 ans plus tard. L'auteur de cet expos\'e recommande la lecture de \cite{DaD} sur les travaux de S.Germain et les tatonnements de l'\'epoque pour d\'eveloper une th\'eorie de l'elasticit\'e non-lin\'eaire bi-dimensionnelle. } seront contest\'es par certains de ses contemporain comme S.Poisson. Quelques d\'ecennies plus tard, en 1850, G.Kirkchhoff donne une place rigoureuse au {\it lagrangien} $W$ dans la th\'eorie de l'elasticit\'e moderne comme \'etant l'\'energie libre dune membrane bidimensionnelle. La propri\'et\'e mixte et assez universelle d'\^etre \`a la fois invariant conforme et d'inclure toutes les surfaces minimales dans ses points critiques a fait que $W$ ne cesse de r\'eappara\^itre dans de nombreux domaines des sciences. En dehors de la g\'eom\'etrie conforme ou de l'elasticit\'e non-lin\'eaire on pourrait aussi citer la biologie c\'ellulaire - avec l'\'energie dite d'Helfrich des membranes lipidiques \`a deux couches \cite{He}  -
ou la relativit\'e g\'en\'erale - c'est le terme principal de la {\it masse de Hawking}....etc.

L'identit\'e (\ref{invariance}) se g\'en\'eralise en fait \`a toute vari\'et\'e. On a le th\'eor\`eme suivant (d\'emontr\'e dans un cadre g\'en\'eral dans \cite{Ch}).

\begin{prop}
\label{inv-conf}
Soit $\Sigma^2$ une vari\'et\'e orient\'e ferm\'ee bi-dimensionnelle et soit $\vec{\Phi}$ une immersion 
de $\Sigma^2$ dans une vari\'et\'e riemannienne orient\'ee $(M^m,g)$. Soit $\Psi$ un diffeomorphisme
positif conforme de $(M^m,g)$ dans une autre vari\'et\'e orient\'e $(N^m,k)$ alors nous avons
\be
\label{VI.34}
 W(\vec{\Phi})+\int_{\Sigma^2}\ov{K}^g\ dvol_{\vec{\Phi}^\ast g}=W(\Psi\circ\vec{\Phi})+\int_{\Sigma^2}\ov{K}^k\ dvol_{(\Psi\circ\vec{\Phi})^{\ast} g}\quad.
\ee
 o\`u $\ov{K}^g$ (resp. $\ov{K}^k$) est la courbure sectionnelle du plan tangent $\vec{\Phi}_\ast T\Sigma^2$ dans  $(M^m,g)$  (resp. du plan tangent $\Psi_\ast\vec{\Phi}_\ast T\Sigma^2$ dans $(N^m,k)$).
\end{prop}
Un cas int\'erressant de l'identit\'e pr\'ec\'edente est celui de l'inverse de la projection st\'er\'eographique\footnote{On rappelle que $\pi(x_1,x_2,x_3,x_4)=(1-x_4)^{-1}(x_1,x_2,x_3)$ est une transformation conforme.} $\Psi=\pi^{-1}$ de ${\R}^3$ dans la sph\`ere unit\'e tridimensionnelle $S^3$ de ${\R}^4$  : soit $\vec{\Phi}$ une immersion
d'une surface ferm\'ee dans ${\R}^3$ alors
\be
\label{s3-r3}
W(\vec{\Phi})=\int_\Sigma |\vec{H}_{\vec{\Phi}}|^2\ dvol_{g_{\vec{\Phi}}}=\int_\Sigma [ |\vec{H}_{\Psi\circ\vec{\Phi}}|^2+1]\ dvol_{g_{\Psi\circ\vec{\Phi}}}=W(\Psi\circ\vec{\Phi}) + \mbox{Aire}(\Psi\circ\vec{\Phi})
\ee
Pour une immersion $\vec{\Phi}$ \`a valeur dans $S^3$ on d\'efinit alors l'\'energie dite de {\it Willmore} - voir plus bas - par
\[
{\mathcal W}(\vec{\Phi}):=\int_\Sigma [ |\vec{H}_{\vec{\Phi}}|^2+1]\ dvol_{g_{\vec{\Phi}}}
\]
Grace \`a la proposition~\ref{inv-conf},  nous savons que ce lagrangien est invariant par composition par des diff\'eomorphismes conformes de $S^3$.
Les surfaces {\it minimales ferm\'ees} de $S^3$, dont on sait qu'elles sont nombreuses\footnote{Les immersions minimales de type {\it Alexandrov} se construisent ais\'emment (\cite{Bre}). Voir sinon les travaux de H.B.Lawson \cite{La}, de H. Karcher, U. Pinkall et I. Sterling, \cite{KPS}, de J. Choe, M. Soret \cite{CS}, ou de N. Kapouleas and D. Wiygul \cite{KW} sur l'existence de plongements minimaux dans $S^3$.}, sont des minima absolus de $W$ et des points critiques de l'aire par d\'efinition. On d\'eduit de l'identit\'e pr\'ec\'edente que les projections st\'er\'eographiques des {\it surfaces minimales ferm\'ees} de $S^3$ sont des {\it surfaces minimales conformes ferm\'ees} de ${\R}^3$ sans \^etre des {\it surfaces minimales} de ${\R}^3$ car il n'en existe pas.

Comme il n'existe pas de {\it surfaces minimales ferm\'ees} dans ${\R}^m$ on imagine facilement que le minimum absolu de $W$ pour toute {\it surface ferm\'ee} ne peut \^etre $0$. Pour toute immersion $\vec{\Phi}$ d'une surface ferm\'ee orient\'e $\Sigma$ L.Simon \'etablit la formule de monotonie suivante (voir \cite{Sim}), pour tout $x^0\in{\R}^m$, et tout $0<t<T<+\infty$
\[
\begin{array}{l}
\ds T^{-2}\,\mbox{Aire}(\vec{\Phi}^{-1}( B_T(\vec{x}^0)))-t^{-2}\,\mbox{Aire}(\vec{\Phi}^{-1}(B_t(\vec{x}^0)))\\[5mm]
\ds=\int_{M\cap B_T(\vec{x}^0)\setminus B_t(\vec{x}^0)}\ \lf|\frac{\pi_{\vec{n}}(\vec{x}-\vec{x}^0)}{|\vec{x}-\vec{x}^0|^2}+\frac{\vec{H}}{2}\rg|^2\ dvol_g
\ds-\frac{1}{4}\int_{\vec{\Phi}^{-1}(B_T(\vec{x}^0)\setminus B_t(\vec{x}^0))}\ |\vec{H}|^2\ dvol_g\\[5mm]
\ds-\frac{1}{T^2}\int_{ \vec{\Phi}^{-1}(B_T(\vec{x}^0))}\ <\vec{x}-\vec{x}^0,\vec{H}>\ dvol_g+\frac{1}{t^2}\int_{\vec{\Phi}^{-1}(B_t(\vec{x}^0))}\ <\vec{x}-\vec{x}^0,\vec{H}>\ dvol_g
\end{array}
\]
En faisant tendre $t$ vers 0 et $T$ vers l'infini on d\'eduit ais\'ement l'in\'egalit\'e dite de Li et Yau
\[
\forall\, x^0\in {\R}^m\quad\quad \theta(|\Phi(\Sigma)|,x^0)=\lim_{t\rightarrow 0}\frac{\mbox{Aire}(\vec{\Phi}^{-1}(B_t(\vec{x}^0)))}{\pi t^2}\le \frac{W(\Phi)}{4\pi}
\]
o\`u $\theta(|\Phi(\Sigma)|,x^0)$ est la densit\'e au point $x^0$ du {\it varifold } donn\'e par l'image de $\Sigma$ par $\vec{\Phi}$. Dans le cas pr\'esent o\`u $\vec{\Phi}$ est suppos\'e \^etre
une immersion r\'eguli\`ere, $\theta(|\Phi(\Sigma)|,x^0)$ est exactement le nombre d'images r\'eciproques de $x^0$ par $\vec{\Phi}$. On a donc le r\'esultat suivant.
\begin{theo}
\label{li-yau}[\cite{LY}]
Soit $\vec{\Phi}$ l'immersion d'une surface ferm\'ee $\Sigma$ dans ${\R}^m$. S'il existe un point $x^0$ de ${\R}^m$ ayant $k$ images r\'eciproques par ${\Phi}$ alors
\[
W(\vec{\Phi})\ge 4\pi k\quad.
\] 
\end{theo}
De ce r\'esultat on d\'eduit ainsi que pour toute immersion d'une surface ferm\'ee on a
\[
W(\vec{\Phi})\ge 4\pi
\]
La minoration $4\pi$ est atteinte par la sph\`ere $S^2$ dans ${\R}^3$ et on d\'emontre sans difficult\'e, toujours \`a partir de la formule de monotonie, que seule la sph\`ere et ses images
par des translations, des rotations et des dilatation satisfont cette borne inf\'erieure. De ce th\'eor\`eme on d\'eduit aussi le corollaire suivant

\begin{coro}\cite{LY}
\label{plongement}
Soit $\vec{\Phi}$ une immersion d'une surface ferm\'ee satisfaisant
\[
W(\vec{\Phi})<8\pi
\]
 alors $\vec{\Phi}$ est un {\it plongement}.
 \end{coro}

Cette minoration de $W$ par $4\pi$ pour toute surface ferm\'ee rappelle dans un certain sens la minorarion par $2\pi$ de l'int\'egrale de la courbure d'une courbe ferm\'ee dans l'espace euclidien avec \'egalit\'e
si et seulement si la courbe est plane et convexe. L'int\'egrale de la courbure d'une courbe ferm\'ee est invariante par dilatation et est de ce fait le lagrangien correspondant  \`a $W$ pour les courbes.
Un r\'esultat bien connu de J.Milnor affirme que  l'int\'egrale de la courbure d'une courbe ferm\'ee et nou\'ee dans un espace euclidien est sup\'erieure \`a $4\pi$. Par analogie, et au vue aussi du th\'eor\`eme~\ref{li-yau} on peut imaginer que des hypoth\`eses suppl\'ementaires sur la compl\'exit\'e de la topologie de $\Sigma$ ou sur la compl\'exit\'e de l'immersion $\vec{\Phi}$ - classe conforme, classe d'isotopie...- devrait donner des minorations plus \'elev\'ees que $4\pi$.

Dans les d\'ecennies qui ont suivi le travail de W.Blaschke presqu'aucun r\'esultat a concern\'e le lagrangien $W$ et ses points critiques, les {\it surfaces minimales conformes}. La raison probablement \'etait que
tr\`es peu de {\it surfaces minimales conformes} qui ne soient pas juste des transformations conformes de {\it surfaces minimales} n'\'etaient connus \`a l'\'epoque et on ne savait pas comment en produire de nouvelles. L'article de 1965 de T.J.Willmore, qui semble-t-il ignorait l'existence du travail de Blaschke, a relanc\'e l'\'etude du lagrangien $W$, qui est devenu {\it l'\'energie de Willmore} et de ses points critiques qui sont devenus les {\it surfaces de Willmore}. Dans son article Willmore formule une conjecture qui est l'objet principal de cet expos\'e
\begin{conj}
\label{willmore}\cite{Wi}
Soit $\vec{\Phi}$ une immersion du tore bi-dimensionnel $T^2$  dans ${\R}^3$ alors on a la minoration suivante
\be
\label{will}
W(\vec{\Phi})\ge 2\pi^2
\ee
avec \'egalit\'e si et seulement si $\vec{\Phi}(\Sigma)$ est , modulo l'action des transformations conformes,
l'immersion axiallement sym\'etrique not\'ee  $T_{Wil}$ et obtenue en tournant autour de l'axe $Oz$ le cercle vertical contenu dans le plan $Oxz$ de centre $(\sqrt{2},0,0)$ et de rayon 1. 
\end{conj}
 Pour \'etayer sa conjecture, T.J.Willmore d\'emontre que le tore $T_{Wil}$, qui sera appel\'e {\it tore de Willmore} apr\`es son travail, est bien un point critique stable de $W$ et que l'in\'egalit\'e (\ref{will}) est bien 
 vraie pour toutes les immersions axiallement sym\'etriques.
 
 La {\it conjecture de Willmore} a stimul\'e de nombreuses contributions \`a la th\'eorie des surfaces qui vont au-del\`a de la conjecture elle m\^eme. Il serait difficile de rendre compte dans cet expos\'e de toute l'activit\'e qu'elle a g\'en\'er\'e. Nous rappellerons seulement quelques r\'esultats qui nous seront utiles pour pr\'esenter la preuve de Marques et Neves dans la section suivante.
 
 Dans un article qui a \'et\'e pr\'ecurseur dans l'analyse des surfaces de Willmore\footnote{L'analyse des surfaces de Willmore est un domaine actuellement tr\`es dynamique de l'analyse g\'eom\'etrique dont la pr\'esentation n\'ec\'essiterait de sortir du cadre de cet expos\'e malheureusement.} \cite{Sim1}, L.Simon d\'emontre, en utilisant la th\'eorie des {\it varifolds} d'Almgren,  que le minimum de $W$ est bien atteint parmi toutes les immersions de $T^2$. Une approche possible pour d\'emontrer la conjecure serait alors d'\'etudier les niveaux d'\'energie occup\'es par les points critiques de $W$ et d'aprofondir notre compr\'ehension
 des {\it surfaces de Willmore}. Une telle approche s'est av\'er\'ee tr\`es \'efficace dans le cas des {\it sph\`eres de Willmore} de $S^3$. Dans \cite{Br} R.Bryant d\'emontre que l'{\it application de Gauss conforme}\footnote{L'application de Gauss conforme associe en chaque point de l'immersion la sph\`ere bi-dimensionnelle tangente dans $S^3$ et de m\^eme courbure moyenne que la courbure moyenne de l'immersion.} d'une immersion dans $S^3$ est une {\it application harmonique} \`a valeur dans la {\it sph\`ere Lorentzienne} $S^{1,3}\subset {\R}^{1,4}$.
 \`A partir de cette application harmonique qui est aussi conforme, donc minimale, en utilisant la repr\'esentation de Weierstrass, R.Bryant construit explicitement une {\it forme quartique holomorphe}\footnote{C'est \`a dire une {\it section holomorphe} du fibr\'e $(\wedge^{(1,0)}T\Sigma)^{\otimes 4}$} dont la nullit\'e est \'equivalente au fait que la {\it surface de Willmore} dans $S^3$ est issue
 de l'image r\'eciproque d'une {\it surface minimale} non compacte de ${\R}^3$ ayant un nombre fini de bouts plong\'es et plans. On v\'erifie facilement que  l'\'energie de Willmore  un multiple de $4\pi$, le degr\'e de multiplicit\'e correspondant au nombre de bouts de la {\it surface minimale}.  Si $\Sigma$ est de genre nul  toute {\it forme quartique holomorphe} y est nulle et donc on en d\'eduit que l'energie de Willmore d'une sph\`ere de Willmore dans $S^3$ ou ${\R}^3$ est n\'ecessairement un multiple de $4\pi$ - R.Bryant d\'emontre aussi que les niveaux de $4\pi {\N}^\ast$ ne sont pas tous atteints. Cet argument malheureusement ne fonctionne pas pour les tores ou les autres surfaces ferm\'ees de genre sup\'erieur qui poss\`edent des {\it formes quartiques holomorphes} non-triviales. Les tentatives de d\'emontrer la conjecture de Willmore par une repr\'esentation de type Weierstrass \`a partir du travail de Bryant a stimul\'e de nombreux travaux, avec des contributions importantes de F.Burstall, J.Dorfmeister, F.Pedit, U.Pinkall...etc  mais n'a jusqu'\`a pr\'esent pas abouti.

Dans l'article \cite{LY}, que nous avons d\'eja mentionn\'e plus haut, P.Li et S.T.Yau introduisent la notion de {\it volume conforme}. Soit ${\mathcal M}(S^3)$ le {\it groupe de M\"obius} des diff\'eomorphismes conformes de la sph\`ere $S^3$. \'Etant donn\'ee l'immersion $\vec{\Phi}$ d'une surface ferm\'ee $\Sigma$ on d\'efinit
\[
V_c(\vec{\Phi},3)=\sup_{\Psi\in {\mathcal M}(S^3)}\mbox{Aire}(\Psi\circ\vec{\Phi})
\]
la quantit\'e $V_c$ est appel\'ee {\it volume conforme} de l'immersion $\vec{\Phi}$. Comme l'{\it\'energie de Willmore} d'une immersion quelconque de $\Sigma$ dans $S^3$ majore l'aire de cette immersion
et comme l{\it'\'energie de Willmore} est invariante par composition par des diff\'eomorphismes conformes on en d\'eduit la proposition suivante.
\begin{prop}\cite{LY}
\label{vol-conf}
Soit $\Sigma$ une surface orient\'ee ferm\'ee et soit $\vec{\Phi}$ une immersion de $\Sigma$ dans $S^3$ alors
\[
V_c(\vec{\Phi},3)\le {\mathcal W}(\vec{\Phi})
\] 
avec \'egalit\'e si et seulment si $\vec{\Phi}$ est l'image par un diff\'eomorphisme conforme d'une immersion minimale.
\end{prop}
\'Etant donn\'e  une surface de riemann $(\Sigma,h)$, o\`u $h$ d\'esigne la m\'etrique de courbure de Gauss constante et de volume 1 on d\'efinit
\[
V_c((\Sigma,h),3):=\inf\lf\{V_c(\vec{\Phi},3)\ ;\ \vec{\Phi}\mbox{ est une immersion conforme de }(\Sigma,h)\mbox{ dans }S^3\rg\}
\]
Un des r\'esultats principaux de \cite{LY} est le th\'eor\`eme suivant
\begin{theo}\cite{LY}
\label{li-yau-vol-conf}
Soit $(\Sigma,h)$ une surface de riemann ferm\'ee \'equip\'ee de sa m\'etrique de courbure constante et de volume 1 alors
\be
\label{val-prop}
\la_1(\Sigma,h)\le 2\, V_c((\Sigma,h),3)
\ee
o\`u $\la_1(\Sigma,h)$ est la premi\`ere valeur propre non nulle du Laplacien sur $(\Sigma,h)$.
\end{theo}
La connaissance explicite des premi\`ere valeurs propres du laplacien des {\it tores plats} permet, en combinant  la proposition~\ref{vol-conf} et le th\'eor\`eme~\ref{li-yau-vol-conf}, d'obtenir 
une minoration explicite de l'{\it\'energie de Willmore} pour toute immersion de tore dont on conna\^\i t la classe conforme. En particulier, Li et Yau parviennent \`a d'\'emontrer la minoration
(\ref{will}) pour un sous domaine de l'espace des modules. Ce sous-domaine sera \'etendu par S.Montiel et A.Ros qui obtiennent de nouvelles minorations pour le {\it volume conforme} des immersions de tores, \cite{MR}.
N\'eanmoins, apr\`es les publications de ces deux travaux, importants pour la conjecture de Willmore, il restait encore une grande partie de l'espace de modules pour laquelle l'in\'egalit\'ee (\ref{will}) etait encore ouverte.
D'autres contributions comme \cite{Ro} ou \cite{To} permettent de d\'emontrer que la conjecture de Willmore est vraie pour toute surface de $S^3$ invariante par l'application antipodale $x\rightarrow -x$.

L'invariance conforme de ${\mathcal W}$ donne que le maximum de l'aire pour une {\it immersion minimale} est atteint exactement pour cette immersion qui r\'ealise donc son volume conforme. C'est justement le cas pour l'image r\'eciproque du {\it tore de Willmore} par la projection st\'er\'eographique $\pi$. Cette image r\'eciproque $\pi^{-1}(T_{Wil})$ n'est autre - modulo rotation - que le {\it tore de Clifford} $T_{Clif}$ donn\'ee par
\[
T_{Clif}:=\lf\{(\sqrt{2})^{-1}(e^{i\theta},e^{i\phi})\in S^{3}\subset {\R}^4\simeq {\C}^2\ ;\ \mbox{ pour tout }(\theta,\phi)\in [0,2\pi)^2\rg\}
\]
et qui est une surface minimale de $S^3$ et on a  $2\pi^2=W(T_{Wil})={\mathcal W}(T_{Clif})=\mbox{Aire}(T_{Clif})$ . Si donc  la conjecture de Willmore \'etait vrai, le tore de Clifford serait en particulier la surface minimale de $S^3$ d'aire minimale parmi toute les surfaces minimales de genre non nul dans $S^3$. C'est cette derni\`ere affirmation que F.Marques et A.Neves vont s'\'efforcer de d\'emontrer rigoureusement.

\subsection{La preuve de la conjecture de Willmore par F.Marques et A.Neves.}

Dans cette sous-section nous rendons compte de la d\'emonstration des r\'esultats suivants de F.Marques et A.Neves, ce qui constituait l'objectif principal de cet expos\'e.

\begin{theo}\cite{MN}
\label{marques-neves-minimal}
Soit $\Sigma$ une surface ferm\'ee  de genre non nul.  Soit $\vec{\Phi}$ une immersion minimale de cette surface dans $S^3$ alors
\[
\mbox{Aire}(\Sigma)\ge 2\pi^2\quad,
\]
avec \'egalit\'e si et seulement si $\Sigma=T^2$ et $\vec{\Phi}(T^2)$ est le tore de Clifford - modulo l'action des transformations rigides. 
\end{theo}
Comme ${\mathcal W}=$Aire pour les surfaces minimales de $S^3$, ce r\'esultat est la cons\'equence du th\'eor\`eme suivant qui, paradoxalement - en apparence  seulement -, sera  d\'emontr\'e dans un deuxi\`eme
temps dans le travail de F.Marques et A.Neves.
\begin{theo}\cite{MN}
\label{marques-neves-willmore}
Soit $\Sigma$ une surface ferm\'ee  de genre non nul. Pour toute immersion de $\Sigma$ dans $S^3$
on a
\[
{\mathcal W}(\vec{\Phi})\ge 2\pi^2\quad,
\]
avec \'egalit\'e si et seulement si $\Sigma=T^2$ et $\vec{\Phi}(T^2)$ est le tore de Clifford - modulo l'action du groupe de M\"obius des transformations conformes de $S^3$.
\end{theo}

Comme on peut le voir ce th\'eor\`eme apporte une r\'eponse bien plus forte que la conjecture de Willmore~\ref{willmore} elle m\^eme telle qu'elle a \'et\'e formul\'ee \`a l'origine car toutes les surfaces
ferm\'ees  de genre non nul sont prises en compte dans ce r\'esultat et pas seulement les tores.

Avant de pr\'esenter les id\'ees principales de la preuve du th\'eor\`eme~\ref{marques-neves-minimal}, il est l\'egitime de s'interroger  sur la situation des surfaces minimales ferm\'ees de genre nul dans $S^3$. 
En 1966 F.Almgren  d\'emontre le r\'esultat suivant que nous allons utiliser plus bas.
\begin{theo}\cite{Alm1}
\label{alm-minimal-spheres}
Soit $\vec{\Phi}$ une immersion minimale d'une surface ferm\'ee de genre nul dans $S^3$ alors $\vec{\Phi}(S^2)$ est une sph\`ere g\'eod\'esique, isom\'etrique \`a $S^2\subset S^3$.
\end{theo}
La d\'emonstration de ce r\'esultat est relativement \'el\'ementaire. Modulo un changement de param\'etrization, le th\'eor\`eme d'uniformisation nous permet de supposer que $\vec{\Phi}$ est {\it conforme et minimale}  de $S^2$ dans $S^3$. L'{\it identit\'e de Codazzi} s'\'ecrit $\ov{\p}h^0=g_{\C}\otimes \p H$ o\`u $h^0$ est la forme quadratique de {\it Weingarten}\footnote{ Le vecteur $\vec{n}$ est le vecteur de Gauss de l'immersion, le vecteur  tangent \`a $S^3$ et perpendiculaire \`a l'immersion unit\'e et positif pour l'orientation choisie.} donn\'ee par ${h}^0:= \vec{n}\cdot\p^2_{z^2}\vec{\Phi}\ dz^2$  et $g_{\C}:=|\p\vec{\Phi}|^2 d\ov{z}\otimes dz$. La nullit\'e de $H$ donne que la forme de {\it Weingarten} est holomorphe. Toute forme quadratique holomorphe sur la sph\`ere est nulle ce qui donne finalement que la forme quadratique de {\it Weingarten} est identiquement nulle. Comme la courbure moyenne est nulle elle aussi on d\'eduit que la seconde forme fondamentale de l'immersion $\vec{\Phi}$ est identiquement nulle et donc que $\vec{\Phi}(S^2)$ est totallement g\'eod\'esique.

\medskip

Le point de d\'epart de la d\'emonstration de la conjecture de Willmore par F.Marques et A.Neves est la caract\'erisation suivante du {\it tore de Clifford} parmi toutes les immersions minimales due \`a F.Urbano.
\begin{theo}\cite{Ur}
\label{urbano}
L'indice de Morse de l'immersion minimale dans $S^3$ d'une surface ferm\'ee de genre non nul est sup\'erieur ou \'egal \`a 5. Il est exactement \'egal \`a $5$ si est seulement si cette immersion
est l'image par une isom\'etrie du tore de Clifford.
\end{theo}

L{\it'indice de Morse} est \'egal \`a la dimension de l'\'espace des vecteurs propres de la d\'eriv\'ee seconde de l'aire, l'op\'erateur de Jacobi de la surface, de valeur propre n\'egative. 
Tout d'abord le fait que l'indice de Morse ne soit pas nul n'est pas une compl\`ete surprise. En effet nous avons vu plus haut que toute {\it immersion minimale} r\'ealise le maximum de son {\it volume conforme}, donc l'action du groupe conforme ne  peut que donner des direction infinit\'esimale n\'egative.

Le fait que toute surface minimale non totallement g\'eod\'esique ait un indice de Morse au moins \'egal \`a 5 se comprend aussi sans trop d'\'efforts. Un calcul tr\`es classique en th\'eorie des surfaces minimales (voir par exemple \cite{CM}) donne que pour toute immersion minimale $\vec{\Phi}$ d'une surface ferm\'ee et 
pour toute perturbation $\vec{w}=w\,\vec{n}$ on a
\be
\label{jacobi}
\lf.\frac{d^2 \mbox{Aire}(\vec{\Phi}+t\,\vec{w})}{dt^2}\rg|_{t=0}=\,\int_\Sigma\lf[|dw|^2_{g_{\vec{\Phi}}}-(|\vec{\mathbb I}_{\vec{\Phi}}|^2+2)\ |w|^2\rg]\,dvol_{g_{\vec{\Phi}}}\quad,
\ee
L'op\'erateur de Jacobi\footnote{$\Delta_{g_{\vec{\Phi}}}$ d\'esigne l'{\it op\'erateur de Laplace Beltrami} - positif - associ\'e \`a la m\'etrique $g_{\vec{\Phi}}$ donn\'e en coordon\'ees locale par
\[
\Delta_{g_{\vec{\Phi}}}w:=-(det (g^{kl}))^{1/2}\sum_{i,j=1}^2\p_{x_i}\lf[det(g_{kl})^{1/2} g^{ij}\p_{x_j}w\rg]
\]
o\`u on rapelle que $g_{\vec{\Phi}}=g_{ij}\ dx_i\otimes dx_j$ et $(g^{kl})$ est l'inverse de la matrice $(g_{kl})$. } est donc ${\mathcal L}_{\vec{\Phi}}w=\Delta_{g_{\vec{\Phi}}}w-(|\vec{\mathbb I}_{\vec{\Phi}}|^2+2)w$. 
Pour tout $v$ dans $B^4$ on introduit la transformation conforme de $S^3$ donn\'ee par
\be
\label{moebius}
\forall z\in S^3\quad\quad F_v(z):=(1-|v|^2)\frac{z-v}{|z-v|^2}-v
\ee
On v\'erifie que pour tout $a\in {\R}^4$ on a $d\, F_{ta}(z)/dt|_{t=0}= 2a\cdot z\ z-2\, a$. Donc la variation normale \`a une surface correspondant \`a l'action infinit\'esimale du groupe conforme - modulo l'action des isom\'etries - est donn\'ee par
\[
w_a:=-2\ a\cdot \vec{n}
\]
L'application de Gauss d'une {\it immersion minimale} dans $S^3$ est harmonique\footnote{Une application harmonique $u$ de $(\Sigma,g)$ dans $S^3$ est un point critique de l'\'energie de Dirichlet $\int_\Sigma|du|^2_g dvol_g$ et satisfait $\Delta_gu-u |du|^2_g=0$.} donc pour toute {\it immersion minimale} on a
\[
{\mathcal L}_{\vec{\Phi}} w_a=-2\,a\cdot\lf[ \Delta_{g_{\vec{\Phi}}}\vec{n}-\vec{n}\ [|d\vec{n}|^2_{g_{\vec{\Phi}}} +2] \rg]=-2\, w_a
\]
Les fonctions $w_a$ r\'ealisent ainsi des fonctions propres de l{\it op\'erateur de Jacobi} et comme la surface minimale n'est pas une sph\`ere g\'eod\'esique, ce que l'on sait grace au th\'eor\`eme de F.Almgren~\ref{alm-minimal-spheres}, cet espace de fonctions propres est de dimension 4. La valeur propre $-2$ a donc une multiplicit\'e au moins \'egale \`a 4. Donc $-2$ ne peut \^etre la premi\`ere valeur propre qui est de mltiplicit\'e 1. Ainsi il existe au moins $5$ vecteurs propres ind\'ependant de ${\mathcal L}_{\vec{\Phi}}$ ayant des valeurs propres n\'egatives et donc l'indice de toute surface minimale non-g\'eod\'esique est au moins 5.

 On v\'erifie ais\'ement  que
 \[
 {\mathcal L}_{T_{Clif}}w=\Delta_{T_{Clif}} w-4\, w
 \]
 Donc la premi\`ere  valeur propre de 'op\'erateur de Jacobi du  {\it tore de Clifford} est $\la_1(T_{Clif})=-4$ avec pour fonctions propres les fonctions constantes. On observe que cette premi\`ere valeur propre
 est g\'en\'er\'ee par la variation infinit\'esimale donn\'ee par
 \[
\vec{\Phi}_t\ :\  (\theta,\phi)\longrightarrow \ (\cos \,(t/4+\pi/4)\  e^{i\theta}, \sin\,(t/4+\pi/4)\  e^{i\phi})\quad\forall \,t\ \in (-\pi,\pi).
 \]
Les immersions $\vec{\Phi}_t(T^2)$ sont les {\it tores \`a courbure moyenne constante} obtenus en relevant les diff\'erents cercles horizontaux de $S^2$ par la fibration de Hopf\footnote{La fibration de Hopf
associe \`a $(z_1,z_2)\in {\C}^2\simeq{\R}^4$ la droite complexe $[z_1,z_2]\in {\C}{\mathbb P}^1\simeq S^2$ passant par ce point.} Le {\it tore de Clifford} \'etant le relev\'e de l'\'equateur.
La famille $\vec{\Phi}_t(T^2)$ correspond au feuilletage de $S^3\setminus\{\mbox{p\^oles nords et sud}\}$ correspondant  \`a
\[
\Sigma_t:=\p\{x\in S^3\ ;\ d(x)<t\}
\]
o\`u $d(x)=\pm \mbox{dist}_{S^3}(x,T_{Clif})$ est la distance sign\'ee \`a $T_{Clif}$  pour la m\'etrique standard sur $S^3$ que l'on note $\mbox{dist}_{S^3}$ et $d$ est compt\'ee soit positivement soit n\'egativement suivant que l'on soit dans une des deux composantes connexes de $S^3\setminus T_{Clif}$. Il est naturel de penser que. modulo
les d\'eformations li\'ees \`a l'action des diff\'eomorphismes conformes et des d\'eformations $\Sigma_t$, le {\it tore de Clifford} minimise l'aire parmi toutes les surfaces de genre non nul.
F.Marques et  A.Neves imaginent donc une proc\'edure de {\it min-max } qui puisse rendre rigoureuse une telle affirmation.

Un glissement de notation est d\'esormais utile afin d'\^etre compatible avec \cite{MN} : $\Sigma$ ne va plus noter une surface abstraite mais une immersion de celle ci.

Pour toute surface {\it plong\'ee}\footnote{Pour d\'emontrer les th\'eor\`emes~\ref{marques-neves-minimal} et ~\ref{marques-neves-willmore} Il est suffisant de ne consid\'erer que les immersions sous le niveau d'\'energie de Willmore ${\mathcal W}(\Sigma)<8\pi$ qui sont donc n\'ecessairement 
plong\'ees comme l'affirme le corollaire~\ref{plongement}.} $\Sigma$ dans $S^3$, on introduit la {\it famille canonique primitive} 
\[
\Sigma_{(v,t)}:=\p\{x\in S^3\ ;\ d_v(x)<t\}\quad\forall\ v\in B^4\ \forall\ t\in [-\pi,\pi] 
\]
o\`u $d_v$ d\'esigne la fonction distance sign\'ee d\'efinie par $d_v(x)=dist_{S^3}(x,F_v(\Sigma))$ pour $x\in F_v(A)$ et $d_v(x)=-dist_{S^3}(x,F_v(\Sigma))$ pour $x\in S^3\setminus F_v(A)$ et $A$ est l'ouvert donn\'e par une des composantes connexes de $S^3\setminus \Sigma$. Les courants d'int\'egration sur $\Sigma_{(v,t)}$, que l'on note $[\Sigma(v,t)]$ d\'efinissent des {\it cycles entiers rectifiables} de ${\mathcal Z}_2(S^3)$ et on v\'erifie sans probl\`eme que
l'application $\Sigma(v,t)$ de $B^4\times (-\pi,\pi)$ dans les {\it cycles entiers rectifiables} bi-dimensionels \'equip\'es de la topologie {\it b\'emol}, ${\mathcal Z}_2(S^3,{\mathcal F})$, est \underbar{continue}. Dans l'espoir de d\'emontrer que 
le {\it tore de Clifford} a l'aire minimale parmi toute les autres {\it surfaces minimales} de genre non nul on envisage alors de consid\'erer
des probl\`emes de {\it min-max} sur les surfaces plong\'ees de genre non nul  que l'on \'ecrit formellement \`a ce stade
\be
\label{min-max-canonique}
\inf_{\Sigma(v,t)\,\simeq\,\Sigma_0(v,t)}\ \ \sup_{v\in\, B^4\, , \, t\in [-\pi,\pi]}\mbox{Aire}(\Sigma(v,t))\quad.
\ee
Ce probl\`eme , on l'esp\`ere, devrait produire une surface minimale d'indice au plus $5$ et qui serait d'aire plus basse que toutes les autres, grace au th\'eor\`eme d'Urbano on en d\'eduirait le
 th\'eor\`eme~\ref{marques-neves-minimal}.

 Il y a une bonne nouvelle \`a ce stade qui est apport\'ee par le th\'eor\`eme suivant de A.Ros.
\begin{theo}\cite{Ro}
\label{ros}
Pour toute surface plong\'ee $\Sigma$ dans $S^3$ on a la majoration suivante
\be
\label{canonique-borne}
\forall\, t\in (-\pi,\pi)\quad\forall\, v\,\in\, B^4\quad\quad\mbox{Aire}\,(\Sigma(v,t))\le {\mathcal W}(\Sigma)
\ee
Par ailleurs, si $\Sigma$ n'est pas une sph\`ere g\'eod\'esique et si il existe $t\in (-\pi,\pi)$ et $v\,\in \, B^4$ tels que
\[
\mbox{Aire}\,(\Sigma(v,t))= {\mathcal W}(\Sigma)
\]
alors $t=0$ et $\Sigma_v$ est une surface minimale.
\end{theo}
La bonne nouvelle est que le maximum de l'aire est born\'e sur toute {\it famille canonique primitive}.  Ce qui n'est par contre pas clair du tout \`a ce stade c'est la signification de l'\'equivalence d'homotopie
que l'on a not\'ee $\simeq$ dans (\ref{min-max-canonique}). Contrairement aux probl\`emes de {\it min-max} consid\'er\'es dans la section pr\'ec\'edente, par passage \`a la limite $v\rightarrow \p B^4$, le bord du cube de l'application $\Sigma(v,t)$ n'est pas envoy\'ee sur l'\'el\'ement nul de ${\mathcal Z}_2(S^3)$. On a par exemple
\be
\label{bord+}
\lf\{
\begin{array}{l}
\ds\forall\ t\,\in\,(0,\pi)\quad\forall\ p\in A\quad\quad{\mathcal F}-\lim_{v\rightarrow p}\ [\Sigma(v,t)]=[\p B^3_{t}(-p)]\\[5mm]
\ds\forall\ t\,\in\,(-\pi,0]\quad\forall\ p\in A\quad\quad{\mathcal F}-\lim_{v\rightarrow p}\ [\Sigma(v,t)]=0
\end{array}
\rg.
\ee
o\`u $[\p B^3_t(-p)]$ est le courant d'int\'egration sur le bord le la boule g\'eod\'esique de centre $-p$ et de rayon $t$ orient\'ee par son vecteur unitaire sortant. On peut alors imaginer remplacer des d\'eformations d'homotopies de bord nul par de d\'eformations ${\mathcal F}-$continues et \underbar{relatives} pour un bord contraint \`a \'evoluer dans un sous ensemble de ${\mathcal Z}_2(S^3)$. C'est alors qu'intervient l'inconv\'enient principal de la {\it famille canonique primitive} : on a
\be
\label{bord-}
\lf\{
\begin{array}{l}
\ds\forall\ t\,\in\,[0,\pi)\quad\forall\ p\in S^3\setminus\overline{A}\quad\quad{\mathcal F}-\lim_{v\rightarrow p}[\Sigma(v,t)]=0\\[5mm]
\ds\forall\ t\,\in\,(0,\pi)\quad\forall\ p\in S^3\setminus\overline{A}\quad\quad{\mathcal F}-\lim_{v\rightarrow p}[\Sigma(v,t)]=[\p B^3_{\pi+t}(p)]\\[5mm]
\end{array}
\rg.
\ee
Donc en comparant (\ref{bord+}) et (\ref{bord-}) en en consid\'erant des  $p$ respectivement dans $A$ et $S^3\setminus \ov{A}$ qui convergent vers un point de $\Sigma$ de part et d'autre de la surface  on constate que la famile canonique primitive \underbar{ne s'\'etend pas de fa\c con continue} pour la {\it topologie b\'emol} au bord du cube $\p(B^4\times(-\pi,\pi))$.

C'est l\`a qu'intervient une des contributions les plus remarquables du travail de F.Marques et A.Neves : c'est la modification de la {\it famille canonique primitive} en une {\it famille canonique ${\mathcal F}-$continue} jusqu'au bord qui pr\'eserve les propri\'et\'es \'essentielles de la famille primitive. Pr\'ecisemment nous avons le r\'esultat suivant.
\begin{theo}\cite{MN}
\label{famille-canonique}
{\bf[existence de familles canoniques ${\mathcal F}-$continues]}
Soit $\Sigma$ une surface ferm\'ee  de genre non nul  plong\'ee dans $S^3$. Alors il existe une application $\Psi_\Sigma$ ${\mathcal F}-$continue jusqu'au bord de $[0,1]^5$ dans ${\mathcal Z}_2(S^3)$ satisfaisant
les conditions suivantes
\begin{itemize}
\item[i)]
\[
\forall\, x'\in[0,1]^4\quad\Psi_\Sigma(x',0)=\Psi_\Sigma(x',1)=0
\]
\item[ii)] Pour tout $x'\in \p([0,1]^4)$ l'application $t\rightarrow \Psi_\Sigma(x',t)$ r\'ealise un feuilletage de $S^3$ par des sph\`eres  centr\'ees en un point $Q_\Sigma(x')\in S^3$.
L'application $Q_\Sigma$ est appel\'ee ''application centre''. En particulier
\[
\sup_{x\in \p ([0,1]^5)}{\mathbf M}(\Psi_\Sigma(x))= 4\pi\quad.
\]
\item[iii)] Pour tout $x'$ dans $\p([0,1]^4)$ on a $$\Psi_\Sigma(x',1/2)=[\p B^3_{\pi/2}(Q_\Sigma(x'))]$$
et $t=1/2$ est le seul $t$ pour lequel c'est vrai. En particulier on a
\[
\forall\ep>0\quad\exists \delta>0\quad{\mathbf F}(|\Psi_\Sigma(x,t)|,{\mathcal T})<\ep\ \Rightarrow\ |t-1/2|<\delta\quad.
\]
o\`u ${\mathcal T}$ d\'esigne l'espace des {\it varifolds} g\'en\'er\'es par les grands cercles de $S^3$.
\item[iv)] La restriction de l'application
\[
\Psi_\Sigma\ :\ \p([0,1]^4)\times [0,1]\,\longrightarrow\, {\mathcal Z}_2(S^3)
\]
est continue pour la topologie ${\mathbf F}$.
\item[v)] 
\[
[\Sigma]\in \Psi_\Sigma([0,1]^5)
\]
\item[vi)]
\[
\sup_{x\in [0,1]^5}{\mathbf M}(\Psi_\Sigma(x))\le {\mathcal W}(\Sigma)\quad.
\]
\item[vii)]
\[
\lim_{r\rightarrow 0}\sup\lf\{{\mathbf M}(\Psi_\Sigma(x)\res B^3_r(p))\ ;\  x\in[0,1]^5\ p\in S^3\ r\in[0,\pi]\rg\}=0
\]
\end{itemize}
\end{theo}

La construction de la {\it famille canonique ${\mathcal F}-$continue} $\Psi_\Sigma$ \`a partir de la {\it famille canonique primitive} est r\'ealis\'ee comme suit. On part de la {\it famille canonique primitive} qui est bien ${\mathcal F}-$continue de $\ov{B^4}\setminus\Sigma\times [-\pi,\pi]$. On observe
que lorsque $v\in B^4$ converge vers un point du bord de la boule contenu dans la surface $p\in\Sigma$ il existe bien une ${\mathcal F}-$limite pour toute suite contenue dans le plan bi-dimensionel
g\'en\'er\'e par les deux vecteurs orthogonaux $p$ et le vecteur de Gauss de la surface $\vec{n}(p)$ et de la forme
\[
v_n=(1-\sigma_n)\,(\cos \tau_n\ p+\sin \tau_n\ \vec{n}(p))\quad\mbox{ telle que }\frac{\tau_n}{\sigma_n}\mbox{ converge }
\]
cette limite est donn\'ee par
\[
\lim_{n\rightarrow +\infty}[\Sigma(v_n,t)]=[\p B_{\frac{\pi}{2}-\theta+t}(-\sin\theta\,p-\cos\theta\,\vec{n}(p))]\quad\mbox{ o\`u }\quad\theta=\lim_{n\rightarrow +\infty}\arctan\frac{\tau_n}{\sigma_n} 
\]
On observe que la famille $[\p B_{\frac{\pi}{2}-\theta+t}(-\sin\theta\,p-\cos\theta\,\vec{n}(p))]$ des limites de $[\Sigma(v_n,t)]$ r\'ealise un feuilletage continu par rapport \`a $\theta$ par des sph\`eres bi-dimensionnelles qui \`a $t$ fix\'e va bien interpoler la ${\mathcal F}-$limite de $[\Sigma(v,t)]$ pour $v\in S^3\setminus \Sigma$ donn\'ees \`a gauche et \`a droite respectivement par (\ref{bord+}) et (\ref{bord-}).
L'id\'ee est alors la m\^eme que celle de la {\it proc\'edure de blow-up} en g\'eom\'etrie complexe en ''rajoutant'' continuement ces feuilletages de sph\`eres interm\'ediaires en chaque point de $\Sigma$.
On obtient ainsi la famille $\Psi_\Sigma$ pour laquelle on v\'erifie i)...vii).

Pour toute surface plong\'ee $\Sigma$ on consid\`ere alors la classe d'homotopie relative $\pi_\Sigma$ des applications $\Xi$ qui sont ${\mathcal F}-$continues de $[0,1]^5$ dans ${\mathcal Z}_2(S^3)$ se d\'eformant
${\mathcal F}-$continuement tout en maintenant la d\'eformation fixe sur le bord $\p[0,1]^5$ et on introduit le probl\`eme de {\it min-max} suivant
\be
\label{min-max-relatif}
{\mathbf L}(\pi_\Sigma):=\inf_{\Xi\in \Pi_\Sigma}\ \max_{ x\in[0,1]^5}{\mathbf M}(\Xi(x))
\ee
La premi\`ere question qui se pose alors est : le probl\`eme de {\it min-max} n'est-il pas trivial\footnote{Cette condition de non-trivialit\'e du {\it min-max} dans le cas d'homotopies relatives correspond \'evidemment \`a la condition de non-trivialit\'e du {\it min-max} dans le cas du travail de 
Pitts pr\'esent\'e plus haut et qui se caract\'rise cette fois par ${\mathbf L}(\Pi)>0$.}
 ?  c'est \`a dire a-t-on bien
\be
\label{min-max-triv}
{\mathbf L}(\pi_\Sigma)>\max_{x\in \p [0,1]^5}{\mathbf M}(\Psi_\Sigma(x))=4\pi\quad\quad ?
\ee
Dans la famille canonique ainsi construite se ''cache'' une rigidit\'e topologique  qui sera \'essentielle pour assurer que le probl\`eme de {\it min-max} n'est pas trivial. 
\begin{theo}\cite{MN}
\label{degre} L'application ''centre'' d\'efinie dans le th\'eor\`eme~\ref{famille-canonique}, vue comme application de $\p[0,1]^4\simeq S^3$ dans $S^3$ a un degr\'e \'egale au genre de la surface
\[
\mbox{deg}\,(Q_\Sigma)=\mbox{genre}\,(\Sigma)\quad.
\]
\end{theo}
Ce r\'esultat se d\'emontre par un calcul relativement  \'el\'ementaire - on connait quasiment explicitement $Q_\Sigma$ en fonction de la surface $\Sigma$ - et fait appel essentiellement au th\'eor\`eme de Gauss Bonnet. C'est cependant est une propri\'et\'e cruciale des {\it familles canoniques ${\mathcal F}-$continues} d\'ecouverte par F.Marques et A.Neves qui va assurer la non-trivialit\'e du probl\`eme de {\it min-max} (\ref{min-max-triv}). On a pr\'ecisemment le th\'eor\`eme suivant.
\begin{theo}\cite{MN}
\label{non-trivialite}
Soit $\Sigma$ une surface immerg\'ee ferm\'ee de genre non nul alors
\[
{\mathbf L}(\pi_\Sigma)>4\pi\quad.
\]
\end{theo}
La preuve de ce th\'eor\`eme dans \cite{MN} consiste \`a rendre rigoureux le  raisonnement tr\`es approximatif suivant. Supposons que la {\it largeur du min-max}, ${\mathbf L}(\pi_\Sigma)$, soit \'egale \`a $4\pi$. On consid\`ere alors une suite $\Xi_n$ homotope \`a la {\it famille canonique} $\Psi_\Sigma$ telle que
\[
4\pi<\max_{x\in  [0,1]^5}{\mathbf M}(\Xi_n(x))\le 4\pi+1/n
\]
Pour tout $v\in [0,1]^4$ la famille $t\rightarrow\Xi_n(v,t)$ r\'ealise un balayage de $S^3$ c'est \`a dire elle est homologue au g\'en\'erateur de $H^3(S^3)$. Il existe donc un $t=t_v$ tel que $\Xi_n(v,t_v)$ est le bord
d'un sous ensemble de $S^3$ de volume \'egal \`a  la moiti\'e de $\mbox{vol}(S^3)$. Lorsque $v\in \p[0,1]^4$ nous savons par exemple que $t_v=1/2$. Comme $\Xi_n(v,t_v)$ est d'aire presque \'egale \`a $4\pi$, elle r\'ealise donc presque l'optimum isop\'erim\'etrique et $|\Xi_n(v,t_v)|$ est de ce fait tr\`es proche pour la topologie ${\mathbf F}$ d'une sph\`ere g\'eod\'esique dont l'espace est not\'e ${\mathcal T}$. On construit ainsi  une 
projection ${\mathbf F}-$continue de l'application $v\in [0,1]^4\rightarrow|\Xi_n(v,t_v)|$ dans ${\mathcal T}\simeq{\R}{\mathbb P}^3$. Une telle projection r\'ealise une extension ${\mathbf F}-$continue
de l'application $x'\in \p[0,1]^4\,\rightarrow|\Psi_\Sigma(x',1/2)|=|[\p B^3_{\pi/2}(Q_\Sigma(x'))]|$ dans ${\mathcal T}\simeq{\R}{\mathbb P}^3$  et $|\Psi_\Sigma(\p[0,1]^4,{1}/{2})|$ serait alors un {\it bord} de ${\R}{\mathbb P}^3$. Ce dernier point  contredirait le fait que le degr\'e de $Q$ est non nul. En effet on v\'erifie aisemment que
\[
\lf|\Psi_\Sigma\lf(\p[0,1]^4,1/2\rg)\rg|=2\, \mbox{deg}\,(Q_\Sigma)\ |[{\R}{\mathbb P}^3]|\in H_3({\R}{\mathbb P}^3)\quad.
\]
 Donc, si genre$(\Sigma)\ne 0$ le th\'eor\`eme~\ref{degre} nous dit que  $\mbox{deg}\,(Q_\Sigma)\ne 0$ et donc $\lf|\Psi_\Sigma\lf(\p[0,1]^4,1/2\rg)\rg|\ne 0$ dans $H_3({\R}{\mathbb P}^3)$
 ce qui contredit l'affirmation plus haut que $\lf|\Psi_\Sigma\lf(\p[0,1]^4,1/2\rg)\rg|$ est un bord de ${\R}{\mathbb P}^3$ . On a donc bien $L(\pi_\Sigma)>4\pi$.

\medskip


\noindent{\bf La discretisation ''\`a la Pitts'' du probl\`eme de min-max.}

Cette \'etape est n\'ec\'esaire pour pouvoir adapter la m\'ethode de Pitts et garantir que le probl\`eme de {\it min-max} est bien atteint. 

Pour toute surface ferm\'ee et plong\'ee $\Sigma$, \`a partir de la {\it famille canonique ${\mathcal F}-$continue} on g\'en\`ere une {\it $(5,{\mathbf M})$-suite d'homotopie canonique} $\varphi_\Sigma=(\varphi_\Sigma^j)_{j\in{\N}}$ de qui ${\mathcal F}-$converge vers $\Psi_\Sigma$, \'egale \`a la restriction  $\Psi^0_\Sigma$ de $\Psi_\Sigma$ sur le bord $\p I(5,j)_0$ et satisfaisant
\be
\label{max}
{\mathbf L}(\varphi)\le \max_{x\in  [0,1]^5}{\mathbf M}(\Psi_\Sigma(x))
\ee
 La possibilit\'e d'une telle discr\'etisation et le passge de la ${\mathcal F}-$continuit\'e \`a la ${\mathbf M}$-continuit\'e discr\`ete avec {\it finesse} tendant vers zero est \`a nouveau illustr\'e
 par la figure.
C'est une \'etape technique et d\'elicate du travail de Marques et Neves que nous ne pouvons d\'ecrire ici. Elle fait usage de nombreuses propri\'et\'es de la {\it famille canonique ${\mathcal F}-$continue} $\Psi_\Sigma$ comme sa ${\mathbf F}-$continuit\'e au bord et la non concentration de la masse vii).

De fa\c con tout \`a fait analogue \'a la definition~\ref{n-m-classe-homotopie} dans le cas de bord nul, on d\'efinit la $(5,{\mathbf M})-$classe d'homotopie des $(5,{\mathbf M})$-suites d'homotopie \`a bord fix\'e \'egal \`a $\Psi^0_\Sigma$. L'ensemble des classes d'homotopies associ\'ees
\`a cette relation d'\'equivalence sont not\'ees $\pi_5^\sharp(Z_2(S^3,{\mathbf M}),\Psi^0_\Sigma)$. Soit $\Pi_\Sigma$ la classe de $\varphi_\Sigma$ dans $\pi_5^\sharp(Z_2(S^3,{\mathbf M}),\Psi^0_\Sigma)$. En suivant une version discrete de la preuve du th\'eor\`eme~\ref{min-max-triv} on d\'emontre que si $\Sigma$ est de genre non nul, on a
\[
{\mathbf L}(\Pi_\Sigma)>\max_{x\in \p [0,1]^5}{\mathbf M}(\Psi_\Sigma(x))=4\pi\quad.
\]
On peut alors appliquer la version ''homotopie relative'' pour le bord fix\'e \'egal \`a $\Psi^0_\Sigma$ du th\'eor\`eme~\ref{Pitt-exist-critical-codim-1} de Pitts et on obtient comme corollaire de ce th\'eor\`eme
et de la construction de la  {\it $(5,{\mathbf M})$-suite d'homotopie canonique} par Marques et Neves.
\begin{coro}
\label{existence}
Soit $\Sigma$ une surface ferm\'ee plong\'ee dans $S^3$ et satisfaisant ${\mathcal W}(\Sigma)<8\pi$ alors il existe une surface minimale plong\'ee $\ti{\Sigma}$ de genre non nul  telle que
\[
{\mathbf L}(\Pi_\Sigma)=\mbox{Aire}(\ti{\Sigma})
\] 
\end{coro}
 Il est tr\`es naturel de penser \`a ce stade que l'indice de $\ti{\Sigma}$ est \'egal \`a la dimension $5$ du domaine des applications consid\'er\'ees. C'est en effet le cas.
 \begin{prop}\cite{MN}
 \label{indice}
 Soit $\ti{\Sigma}$ une surface minimale plong\'ee de $S^3$ telle que
 \[
{\mathbf L}(\Pi_\Sigma)=\mbox{Aire}(\ti{\Sigma})
\]
pour une surface $\Sigma$ de genre non nul satisfaisant ${\mathcal W}(\Sigma)<8\pi$. Alors l'indice de Morse de l'op\'erateur de Jacobi ${\mathcal L}_{\ti{\Sigma}}$ de $\ti{\Sigma}$ est \'egal \`a 5 et $\ti{\Sigma}$ est l'image
par une isom\'etrie du tore de Clifford.
 \end{prop}
La d\'emonstration de ce r\'esultat dans \cite{MN} est assez naturelle. On d\'emontre tout d'abord sans trop d'efforts que l'{\it op\'erateur de Jacobi} ${\mathcal L}_{\ti{\Sigma}}$ est non d\'eg\'en\'er\'e sur l'espace des directions n\'egatives donn\'ees par $a\cdot\vec{n}$ et la fonction constante $1$. En supposant que ${\mathcal L}_{\ti{\Sigma}}$ soit d'indice sup\'erieur ou \'egal \`a $6$, on utilise cette non d\'eg\'en\'erescence pour construire une application $\ti{\Psi}$ de $[0,1]^5$ dans ${\mathcal Z}_2(S^3)$ en d\'eformant  $\Psi_\Sigma$ au voisinage de $\ti{\Sigma}$ le long de cette 6-\`eme hypoth\'etique direction n\'egative de fa\c con \`a avoir
\[
\sup_{x\in [0,1]^5}\ti{\Psi}(x)<\mbox{Aire}(\ti{\Sigma})\quad.
\]
La d\'eformation est \'effectu\'ee avec soin de fa\c con \`a garantir l'existence d'une  {\it $(5,{\mathbf M})$-suite d'homotopie canonique} $\ti{\varphi}=(\ti{\varphi}^j)_{j\in{\N}}$ homotope \`a $\Pi_\Sigma$ et telle que
\[
{\mathbf L}(\ti{\varphi})<\mbox{Aire}(\ti{\Sigma})\quad.
\]
Ceci contredit la minimalit\'e de $\ti{\Sigma}$. Donc $\ti{\Sigma}$ est bien d'indice \'egal \`a $5$ et on peut lui appliquer le th\'eor\`eme d'Urbano~\ref{urbano}.

\medskip

 \noindent{\sc Preuve de la conjecture de Willmore.}

 \medskip
 
 Soit $\Sigma$ une surface plong\'ee de genre non-nul satisfaisant ${\mathcal W}(\Sigma)<8\pi$. On consid\`ere la {\it $(5,{\mathbf M})$-suite d'homotopie canonique} $\varphi_\Sigma=(\varphi_\Sigma^j)_{j\in{\N}}$ associ\'ee
 et sa classe d'homotopie $\Pi_\Sigma$. Les in\'egalit\'es  respectivement vi) du th\'eor\`eme~\ref{famille-canonique} et (\ref{max}) donnent
 \[
 {\mathbf L}(\Pi_\Sigma)\le {\mathbf L}(\varphi_\Sigma)\le{\mathcal W}(\Sigma)\quad.
 \]
 Le corollaire~\ref{existence} nous donne l'existence d'une surface minimale pong\'ee de genre non nul telle que
 \[
 \mbox{Aire}(\ti{\Sigma})= {\mathbf L}(\Pi_\Sigma)\le {\mathbf L}(\varphi_\Sigma)\le{\mathcal W}(\Sigma)
 \]
 La proposition~\ref{indice} nous dit que $\ti{\Sigma}$ est d'indice 5 et donc qu'il s'agit de l'image par une isom\'etrie du {\it tore de Clifford} et donc on a
 \[
 2\pi^2\le{\mathcal W}(\Sigma)\quad.
 \]
Si les in\'egalit\'es pr\'ec\'edentes sont des \'egalit\'e alors, grace \`a (\ref{max}) et la deuxi\`eme partie du th\'eor\`eme~\ref{ros}, nous d\'eduisons que $\Sigma$ est l'image par une transformation conforme du
{\it tore de Clifford}.

\section{Conclusion.}

Il est l\'egitime de se demander ce qu'il en est de la conjecture de Willmore en codimension plus grande que 1 : Qu'en est-il de la borne inf\'erieure de ${\mathcal W}(\Sigma)$ pour des surfaces ferm\'ees de genre
non nul dans $S^m$ pour $m\ge 4$. Les r\'esultats partiels de Li et Yau \cite{LY} ou de Montiel et Ros \cite{MR} s'appliquent dans ce cas et on connait d\'eja un sous ensemble ouvert de l'espace des modules
des tores pour lesquels on a n\'ecessairement ${\mathcal W}(\Sigma)\ge 2\pi^2$. La m\'ethode de {\it min-max} de Marques et Neves semble ne pas s'appliquer en codimension plus grande. En effet, la construction de la {\it famille canonique primitive}, essentielle \`a la construction du probl\`eme de {\it min-max}, est manifestement tr\`es li\'ee\footnote{Ceci est \'evident pour la variable $t$ qui pointe vers la direction la plus n\'egative.} \`a la dimension 3 .
Pitts par ailleurs utilise fortement la codimension 1 afin d'assurer la r\'egularit\'e des surfaces minimales r\'ealisant la {\it largeur}. En codimension plus grande, les surfaces minimales localement minimisantes
peuvent avoir des singularit\'es - les {\it courbes holomorphes} en sont un exemple - et il faudrait utiliser les techniques tr\`es lourdes des r\'esultats d'Almgren-De Lellis-Spadaro pour montrer la r\'egularit\'e en dehors de points isol\'es pour  les surfaces {\it presques minimales}  (ce probl\`eme est encore compl\`etement ouvert).
Il y a aussi un ph\'enom\`ene nouveau en dimension plus grande ou \'egale \`a 4 : le r\'esultat d'unicit\'e d'Almgren n'est plus vrai et il existe par exemple une surface minimale de genre nul
dans $S^4$ qui ne soit pas totalement g\'eod\'esique : la {\it surface de Veronese} en est un exemple. Son aire, \'egale \`a $6\pi$, est comprise entre $4\pi$ et $2\pi^2$ et il faudrait s'assurer qu'un \'eventuel probl\`eme de {\it min max} formul\'e dans le cadre des {\it varifolds} - donc dans un cadre tr\`es faible - ne ''d\'eg\'en\`ere'' pas sur cette surface.


Du c\^ot\'e des m\'ethodes de {\it min-max} mentionnons juste le r\'esultat important obtenu il y a quelques mois par F.Marques et A.Neves en combinant la version ${\Z}_2$ des m\'ethodes d'Almgren-Pitts
et des travaux respectivement de L.Guth \cite{Gu} et M.Gromov \cite{Gr1}, \cite{Gr2} sur les notions de {\it $p-$balayages} ${\Z}_2$ de vari\'et\'e et les estimations de leurs {\it largeurs}.

\begin{theo}
\cite{MN1}
\label{infinite-miniales}
Soit $(M^m,g)$ une var\'et\'e riemannienne compacte de dimension $m$ inf\'erieure ou \'egale \`a $7$ et de courbure de Ricci strictement positive. Alors $M^m$ contient une infinit\'e d'hypersurfaces minimales r\'eguli\`eres, plong\'ees et ferm\'ees.
\end{theo}

\end{document}